\documentclass[aos,MSNbibl,seceqn,dvips]{arximspdf}
\usepackage{mathrsfs}
\usepackage{graphicx}
\doi{10.1214/12-AOS1031} 
\volume{40}
\issue{4}
\pubyear{2012}
\firstpage{2162}
\lastpage{2194}

\makeatletter
\newcommand{\eqref}[1]{(\ref{#1})}
\newtheorem{theorem}{Theorem}
\newtheorem{lemma}{Lemma}
\newtheorem{proposition}{Proposition}
\def\ci{\perp\!\!\!\perp}
\def\F{\mathcal{F}}
\def\PP{\mathscr{P}}
\def\E{\mathcal{E}}
\def\V{\mathcal{V}}
\def\J{\mathcal{J}}
\makeatother

\begin{document}
\begin{frontmatter}

\title{Counterfactual analyses with graphical models based on local independence\thanksref{T1}}
\runtitle{Counterfactual analyses and local independence}
\thankstext{T1}{Supported by The Research Council of Norway, Project:
191460/V50, and The Norwegian Cancer Society, Project 2197685.}

\begin{aug}
\author[A]{\fnms{Kjetil} \snm{R{\o}ysland}\corref{}\ead[label=e1]{kjetil.roysland@medisin.uio.no}}
\runauthor{K. R{\o}ysland}
\affiliation{University of Oslo}
\address[A]{Institute of Basic Medical Sciences\\
Department of Biostatistics\\
University of Oslo\\
Boks 1122 Blindern\\
0317 Oslo\\
Norway\\
\printead{e1}} 
\end{aug}

\received{\smonth{12} \syear{2011}}
\revised{\smonth{7} \syear{2012}}

%
\begin{abstract}
We show that one can perform causal inference in a natural way for
continuous-time
scenarios using tools from stochastic analysis.
This provides new alternatives to the
positivity condition for inverse probability
weighting.
The
probability distribution that would govern the frequency of
observations in
the counterfactual scenario can be characterized in terms of a so-called
martingale problem. The counterfactual and factual probability
distributions may be related through a likelihood ratio given by a
stochastic differential equation. We can perform
inference for counterfactual scenarios based on the original
observations, re-weighted according to this likelihood ratio.
This is possible if the solution of the stochastic differential
equation is uniformly integrable, a property that can be determined
by comparing the corresponding factual and counterfactual short-term
predictions.

Local independence graphs are directed, possibly cyclic, graphs
that represent short-term prediction among sufficiently
autonomous stochastic processes.
We show through an example that these
graphs can be used to identify and provide consistent estimators for
counterfactual parameters in continuous time.
This is analogous to how Judea Pearl uses graphical information to
identify causal effects in finite state Bayesian networks.
\end{abstract}

\begin{keyword}[class=AMS]
\kwd{92D30}
\kwd{62N04}
\kwd{60G44}
\kwd{60G55}
\end{keyword}
\begin{keyword}
\kwd{Causal inference}
\kwd{stochastic analysis}
\kwd{event history analysis}
\kwd{marked point processes}
\kwd{change of probability measures}
\kwd{local independence}
\end{keyword}

\end{frontmatter}

\section{Introduction}\label{sec1}

While randomized controlled trials are the gold standard for
determining the effects of public health interventions or medical
treatments, there are many situations where such trials are unethical,
and it is tempting to turn to registry data or observational studies
for quality assessment of treatments. However, data from such sources
is subject to various selection effects from drop-out due to
underlying health problems to selection of the treatment itself.
These problems have motivated the development of the field of causal
inference, including in particular the area of marginal structural
models \cite{Robins1,royslandmsm} which have seen
applications, for instance, in HIV
cohort studies \cite{Sterne}.
The underlying idea is that observational data can be used to mimic a relevant
hypothetical controlled trial or counterfactual scenario.

In this paper, our primary concern is the possibility of estimating
parameters in a model for the observations from a counterfactual
scenario involving a relevant hypothetical randomized controlled
trial. While the specification of an appropriate model for the
counterfactual observations is an important topic in itself, we will
focus solely on a situation in which such a counterfactual model has
been specified correctly.
It is common to re-weight the observational data in order to
mimic observations coming from the counterfactual scenario.
This is usually referred to as
inverse probability weighting.
Such re-weighting has occasionally been reported
to be too unstable, even inconsistent, for various purposes;
see \cite{HernanCole}.
It is therefore
of great interest to understand when this strategy
actually works. We will provide some rigorous
conditions for such re-weighting to be achievable.
A similar exposition has not been carried out in the
literature before,
except partly in \cite{royslandmsm} and~\cite{HernanCole}.

A
probability distribution on the underlying sample space that would
govern the frequency of observations in
the counterfactual scenario can be characterized in terms of a so-called
martingale problem.
Short-term predictions provide dynamical characterizations of the
various involved modules. A hypothetical direct intervention on a
module would change its dynamics. The nondirectly intervened modules
on the other hand, should have the same dynamical characterization as
in the factual scenario. Martingale problems have been thoroughly
studied in stochastic
analysis; to us one would mean that there would exist well-developed
tools for determining the feasibility of the previous
re-weighting methods. An immediate application of these tools yields,
for instance,
that the probability distribution that would govern the
frequencies of events in the counterfactual situation is
unique if it exists; see Theorem~\ref{uniquenessduringfollow-up} in
the \hyperref[app]{Appendix}.

If the re-weighting is feasible, is it then at all possible
to estimate the parameters of interest in the counterfactual
model from the re-weighted observations?
In other words, are these parameters identifiable?
Pearl's strategy \cite{Pearl} is to take advantage of graphical
structure, in terms of conditional independences, for
identification of causal effects. It was shown in \cite{PearlShpitser,HuangValtorta} and \cite{Sullivan} that this strategy gives
a complete theory in the simpler setting of finite state or
Gaussian--Bayesian networks.
For more complicated settings, this problem is far from
solved. Some results in this direction for time series were given in
\cite{DidelezEichler}. We show
that it is possible to take advantage of local independence graphs for
identification of causal effects in
continuous-time settings. Note,
as this general problem is very hard, we do not
provide a complete theory for
identification of causal effects, only an
example which slightly extends \cite{MartinussenVansteelandt}.

The idea that the counterfactual situation can be assigned
probabilities in a way that is consistent with a purely observational
scheme, is not new. It has also been considered in
the general context of marked
point processes in
[\cite*{didelez3,ArjasParner,Arjas3}, \cite*{parner1999Arjas}] and~\cite{royslandmsm}. We choose a martingale-based
approach, similar to \cite{royslandmsm}.
Note also that graphical models based on local independence and
doubly stochastic Poisson processes were studied thoroughly in
\cite{DidelezRoyalSB}.
Continuous-time counterfactual interventions were also
considered by Lok in \cite{Lok}. She considered
structural nested models in continuous time and applied ideas from structural
equation modeling to survival data. Her strategy differs from
ours in that we take a purely nonparametric point of view, through
change of probability measures.

In Section~\ref{observationalregime} we describe models for
the factual scenario. We then proceed in Section~\ref{actionSection}
with a description of counterfactual variables and distributions. In Section~\ref{counterfactualConstruction},
we give a sufficient condition for
such a counterfactual
distribution to exist, and also a construction based on martingale
methods. In Section~\ref{LI}, we introduce local independence graphs
that play
the same role as directed acyclic graphs usually do in the literature
on causal inference.
In Section~\ref{CE}, we
consider an example where we can identify consistently estimate
controlled direct effects in event history analysis.
Finally, in the \hyperref[app]{Appendix}, we summarize some properties of dual
predictable projections and consider uniqueness of counterfactual
distributions.

\section{The observational regime and autonomous modules}\label{observationalregime}
Eventually, we will consider statistical analyses based on
observations of several
i.i.d. individuals, but first we will consider models for one
``generic'' individual.
We aim to investigate complex systems for each individual
formed by finitely many autonomous modules that develop and influence
each other throughout time. We will not provide a detailed recipe for
building appropriate models,
but simply assume a stochastic model for a generic individual that has
some specific
properties.
\subsection{The underlying probability space and marked point
processes}

We let $\mathcal V$ denote the finite set of modules that form the
system of interest. The possible outcomes of these modules are supposed
to be realized on a probability space $(\Omega, \F, Q)$ with
some additional structure that we will now describe. Note that we do
not assume that the actual frequencies of outcomes
will be governed by the probability measure $Q$. This
measure will only play a role as a ``reference measure.''
The possible ``initial''
outcomes of each module $V$ are given by the outcomes of a corresponding
random variable $V_0$.
The
random variables in this family, which we denote by $\V_0$, are
mutually independent with respect to $Q$.
The intital outcome of each $V \in
\V$ occurs at a, possibly unknown, time point $T(V_0) \leq0$. The
ordering of these time points is assumed to be known.
We moreover let
%
\begin{equation}
p(V_0) := \bigl\{ V'_0 \in
\V_0 | T\bigl(V'_0\bigr) <
T(V_0) \bigr\},
\end{equation}
and sometimes refer to this set as the \textit{past} of
$V_0$.

The outcomes in the follow-up are driven by a multivariate
point process $N$ \cite{JacodMultivariate} on a finite time
interval $[0,T]$. Let $J$ denote the mark space of $N$. This space is
supposed to be Lusin, that is,
a Borel subset in a compact metric space, and equipped
with the Borel $\sigma$-algebra $\J$. We assume that for every module $V$,
there exists a
$J_V \in \J$ such that
%
\begin{equation}
\label{detectableDef} V_t(\omega) = V_0(\omega) + \int
_{J_{V}} \int_0^t h(
\omega,s,x) N(\omega,ds,dx),
\end{equation}
where $h$ is a bounded process on $[0, T]\times J$ that is predictable
with respect to the
filtration generated by $N|_{J_V}$ and $V_0$.
We also assume that $\V_0\ci_Q N$ and that $\coprod_{V \in\V} J_V$
defines a
partition of $J$ such that the restricted point processes
$\{N|_{J_V}\}_{V \in\V}$ are mutually independent with respect to
$Q$.

For each subset $\mathcal W := \{V^1, \ldots, V^d\} \subset\V$, let
$\F_t^{\mathcal W}$ denote the filtration that
is generated by $V_0$ and $N|_{J_{V}}$ for every $V \in\mathcal W$
and also
satisfies the \textit{usual conditions}; see \cite{JacodShiryaev}.
We let $\PP^{ \mathcal W} $ denote the predictable $\sigma$-algebra
generated by $\F_t^{\mathcal W}$ \cite{JacodShiryaev}.
For notational simplicity, we will also write $\F_t^V$ or $\PP^V$
instead of $\F_t^{\{V\}}$ or $\PP^{\{V\}}$, as well as $\F_t$
or $\PP$
instead of $\F_t^{\V}$ or $\PP^{\V}$.

\subsection{The factual distribution}

The actual frequencies of outcomes in the model
are not assumed to be governed by $Q$, but another probability
measure $P$ such that $P \ll Q$
and
%
\begin{equation}
\label{contind} V_0 \ci_P T^{-1}T(V_0)
\setminus\{ V_0\} | p(V_0)
\end{equation}
for every $V_0 \in\V_0$, that is, every $V_0$ is independent
w.r.t. its simultaneous variables, conditionally on the past.
We will refer to property \eqref{contind} as \textit{contemporaneous
independence}; see \cite{DidelezEichler}.
This is useful to us since it provides at least one enumeration $\{
V_0^1, \ldots, V_0^n\} = \V_0 $ such
that
$T(V_0^i) \geq T(V_0^j)$
whenever $i > j$
and
%
\begin{equation}
\label{ordering} E_P \bigl[ f\bigl(V_0^k
\bigr) |V_0^{k-1}, \ldots, V_0^{1}
\bigr] = E_P \bigl[ f\bigl(V_0^k\bigr) | p
\bigl(V_0^k\bigr) \bigr],
\end{equation}
whenever $f$ is a bounded and measurable function and $1 \leq k \leq n$.

The processes in $\V$ are not necessarily mutually independent with
respect to
$P$, but are still sufficiently autonomous for our purpose.
As an immediate manifestation of this autonomy, note that the modules
may not ``switch'' states
simultaneously $P$-a.s.
The reason is that
the processes in $\V$ are associated to disjoint subsets in the
mark space $J$, which cannot occur simultaneously.
We will refer to $P$ as the \textit{factual measure}. Note, however,
as some of the processes in $\V$ may be
latent, the factual measure $P$ is also assumed to
govern the frequency of events that may be unobserved.

\subsection{The factual likelihood ratio and its factorization}

The autonomy imposes a factorization of the
likelihood ratio $\frac{dP }{dQ}$\vspace*{1pt} that will prove to be important
to us. First note that a
repeated use of the Radon--Nikodym theorem
provides a family $\{Z^V_0\}_{V \in\V}$ of nonnegative random
variables such that each $Z^V_0$ is $\F_0^{p(V) \cup\{V\}
}$-measurable and
%
\begin{equation}
E_Q \bigl[ Z^V_0 | \F_0^{p(V)}
\bigr] = 1\quad \mbox{and}\quad \frac{d P |_{\F_0}} {d Q |_{\F_0}} = \prod_{V \in\V}
Z^V_0,\qquad Q\mbox{-a.s.}
\end{equation}

There\vspace*{1pt} is a similar factorization of $\frac{dP }{dQ}$.
Let $U$ denote the dual predictable projection of $N$ with respect
to $Q$ onto the filtration $\F_t$ as in
\cite{JacodMultivariate}.
By Lemma~\ref{predictableprojections} in the \hyperref[app]{Appendix} there exists a
nonnegative and $\PP\otimes\J$-measur\-able
process $\lambda$ such that
\[
E_P \biggl[ \int_J \int_0^T
h( s, x) N(ds, dx) \biggr] = E_P \biggl[ \int_J
\int_0^T h( s, x)\lambda(s, x) U(ds, dx)
\biggr]
\]
for every bounded and $\PP\otimes\J$-measurable
process $h$. As common practice, we mostly omit $\omega$ from equations
in order to be notationally less overwhelming.

We now define
the processes
\[
H^V(t) := 1 + \frac{U( \{t\}, J_V )
- \int_{J_V}\lambda( t , x ) U(\{t\}, dx) } { 1- U( \{t\}, J_V ) }
\]
and
%
\begin{equation}
\label{KV} K^V_t := \int_{J_V} \int
_0^t \lambda(s,x) - H^V(s) \bigl(
N( ds, dx) - U(ds, dx) \bigr).
\end{equation}
By \eqref{strongort2}, we see that
that $\{K^V\}_{V \in\V}$ defines a family of
local $Q$-martingales with respect to the
filtration $\F_t$ such that
%
\begin{equation}
\bigl[ K^V, K^{V'} \bigr] = 0, \qquad Q\mbox{-a.s. for } V \neq
V'.
\end{equation}
The solution of the SDE
%
\begin{equation}
Z_t = Z_0 + \sum_{V \in\V}
\int_0^t Z_{s-}\,dK_s^V
\end{equation}
defines a $Q$-martingale with respect to the filtration $\F_t$
such that
\[
Z_t = \frac{dP|_{\F_t} }{dQ |_{\F_t}}, \qquad Q\mbox{-a.s. }
\]
for every $t \in[0, T]$. This follows directly from \cite{JacodMultivariate}, Theorem~5.1.

We now obtain directly from
Yor's additive formula \cite{Protter}, Theorem II 38, that
%
\begin{equation}
\label{likfactorization} Z_t = \prod_{V \in\V}
Z_t^V,\vadjust{\goodbreak}
\end{equation}
where each $Z^V$ solves an SDE
%
\begin{equation}
\label{ZV} Z_t^V := Z_0^V + \int
_0^t Z_{s-}^V
\,dK_s^V.
\end{equation}

\section{Actions and counterfactual distributions} \label{actionSection}

\def\A{\mathcal{A}}
\def\X{\mathcal{X}}

We assume that we may
directly intervene on a subset of modules $\A\subset\V$ such that
their outcomes are changed. This intervention does not directly affect the
outcomes of the modules in $\X:= \V\setminus\A$. The latter set of
modules will
only be affected indirectly: The conditional distributions of their
short-term behavior, given
the past, will remain the same, while the change of previous
outcomes yields a change of the background these distributions depend on.
We will limit our discussion to actions that are
deterministically dependent on the past. These are sometimes referred
to as conditional actions. Every conditional
action will be represented by a measurable transformation $\theta$ of
the generic state space $( \Omega, \F)$.
We think
of $\theta(\omega)$ as the direct consequence in
the ``counterfactual universe'' where the action $\theta$ was
performed.

Whenever $P'$ is a probability measure on $(\Omega,
\F)$, we let $\theta P'$ denote the push-forward measure over $\theta
$, that is, $\theta P'(F) := P' ( \theta^{-1} ( F) )
$ for every $F \in\F$. Whenever $H$ is an $\F$-measurable random
variable, we let $\theta^* H$ denote the
transformed variable, where $\theta^* H ( \omega) := H( \theta(
\omega))$ for every $\omega\in\Omega$.
We assume that
$\theta$ is ``continuous'' in the sense that the reference measure
$Q$ is quasi-invariant with respect to $\theta$, that is,
%
\begin{equation}
\label{action1} \theta Q \ll Q.
\end{equation}
%
\subsection{Actions and counterfactual distributions at baseline}
Let $V \in\V$ and suppose $\eta$ is an
$\F_0^V$-measurable random variable, and $h$ is a bounded and
$\F_0^{p(V)}$-measurable random variable.
We assume that the outcomes of the not directly intervened part of the
system are left
invariant by the transformation at baseline, that is,
%
\begin{equation}
\label{action2} \theta^* \eta= \eta
\end{equation}
for every $\eta$ and every $V \in\X$.
We furthermore assume that the action depends deterministically on the
past outcomes in the nonintervened system, that is,
whenever $V \in\A$, then
%
\begin{equation}
\label{action4} \theta^* \eta \mbox{ is } \F_0^{p(V)\cap\X}
\mbox{-measurable}
\end{equation}
for every $\eta$.

A probability distribution $P_\theta$ on $(\Omega, \F)$ defines a
\textit{counterfactual distribution at
baseline} if, whenever $V \in\A$, then
%
\begin{equation}
\label{baselineaction2} E_{P_\theta} [ h \eta ] = E_{P_\theta} \bigl[
h \theta^* \eta \bigr],
\end{equation}
and, whenever $V \in\X$, then
%
\begin{equation}
\label{baselineaction1} E_{P_\theta} [ h \eta ] = E_{P_\theta} \bigl[
h \theta^* E_P \bigl[ \eta| \F_0^{p(V)}\bigr]
\bigr]
\end{equation}
for every $\eta$.

Equation \eqref{baselineaction1} means that the short-term behavior
of a directly intervened variable
is simply given by the transformed variable. Its outcome is
deterministically regulated by the past.
Equation \eqref{baselineaction1} means that the conditional
distribution
of
an outcome of a
not directly intervened variable in the counterfactual scenario, given
its past,
coincides with the corresponding
distribution from the factual
scenario.

Note that
Pearl's $\operatorname{do}(X = x)$ may also be interpreted as a
transformation on sample space that
fixes $X$ constantly equal to $x$ and leaves the remaining variables invariant.
This means that our characterization of probability measures
on $(\Omega, \F)$
that
would govern the
frequencies of events in our system if we, \textit{contrary to
the fact}, had applied the hypothetical intervention
strategy, is a reformulation of
Pearl's $do$-operator on Bayesian networks \cite{Pearl}. The present
approach, however, translates more or less directly to
continuous-time settings.

\subsection{Actions and counterfactual distributions in the follow-up
period}

$\!\!\!$When\-ever $Z$ is a
stochastic process on $\Omega$, we let $\theta^* Z$ denote the process
given by the transformed variables $\{\theta^* Z_t\}_{t \in[0,T]}$.
We assume that
$\theta^* N$
defines a marked point process that is adapted to the history
$\{\F_t\}_{t \in[0,T]}$.
The action $\theta$ is thought to force the outcomes
$ N|_{[0, T] \times J_\A }$
into the outcomes of $\theta^* N|_{[0, T] \times J_\A }$, which will
only depend
on the strictly previous behavior of the not directly intervened
system, that is, whenever $B \in J_\A$, then
%
\begin{equation}
\label{action5} \theta^* N_t(B) \qquad\mbox{is predictable w.r.t. } {
\F_t}^\X.
\end{equation}
The outcomes of the not directly intervened part of the system are left
invariant by the transformation during follow-up, that is,
%
\begin{equation}
\label{action3} \theta^* N|_{[0, T] \times J_\X} = N|_{[0, T] \times J_\X}.
\end{equation}

We will say that $P_\theta$ defines a \textit{counterfactual
distribution} if it
defines a counterfactual distribution at baseline, and if whenever
$X$ is process on the form \eqref{detectableDef}
and~$\Lambda$
is an $\F_t$-predictable process of finite variation such that
\[
E_P \biggl[ \int_0^T
h_s \,dX_s \biggr] = E_P \biggl[ \int
_0^T h_s \,d\Lambda_s
\biggr]
\]
for every bounded and $\F_t$-predictable process $h$, then
%
\begin{equation}
\label{follow-upaction1} E_{P_\theta} \biggl[ \int_0^T
h_s \,dX_s \biggr] = E_{P_\theta} \biggl[ \int
_0^T h_s \,d\theta^*
\Lambda_s \biggr] \qquad\mbox {if } V \in\X
\end{equation}
and
%
\begin{equation}
\label{follow-upaction2} E_{P_\theta} \biggl[ \int_0^T
h_s \,dX_s \biggr] = E_{P_\theta} \biggl[ \int
_0^T h_s \,d\theta^*
X_s \biggr] \qquad\mbox{if } V \in\A.
\end{equation}
Note that \eqref{follow-upaction1} means that
$\theta^* \Lambda$ defines the
compensator of $X$ if $V \in\X$, and
\eqref{follow-upaction2} means that
$\theta^* X$ defines the
compensator of $X$, otherwise. This offers an analogous
interpretation as in the baseline setting.
Compensators provide a notion of short-term behavior,
analogously to the previous conditional distributions.
The short-term behavior of a not
directly intervened process in the counterfactual scenario,
based on the past, coincides with the transformed
short-term behavior from the factual
scenario. The short-term behavior of
a directly intervened process is given entirely by the
transformation.

Following \cite{PearlIntroductionToCausalInference}, we will say that
a model consisting of a factual scenario, an action
and a corresponding counterfactual distribution, defines a
causal model
if the counterfactual distribution would
fit the actual corresponding counterfactual scenario.
That $P_\theta$ actually would govern
the frequency of observations for this hypothetical scenario is
generally not testable, and mostly comes down to the
question of \textit{no unmeasured confounding} \cite
{PearlIntroductionToCausalInference}.

\section{Construction of counterfactual distributions} \label
{counterfactualConstruction}

\subsection{Construction at baseline}
We will now construct the counterfactual distribution
in a situation with no follow-up period. The
construction is then closely related to Pearl's framework
\cite{Pearl}. The next result is important and says heuristically that
if the conditional probability, given the past, of observing
outcomes that coincide with counterfactually enforced ones are not
too small, then there exists a counterfactual distribution.
Equation \eqref{baselineactionmeasure} then offers a useful
description of the distribution.
Note that this is a measure
theoretical version of the \textit{truncated factorization formula}
from \cite{Pearl}, (3.10).

\begin{theorem} \label{baselinexistthm}
If
there exists a nonnegative $K \in L^1 (\F_0,P) $ such that
%
\begin{equation}
\label{baselinecondition} \frac{d \theta Q |_{\F_0}}{ d Q |_{\F_0}} \leq K \prod
_{V \in
\A} Z_0^V,\qquad P\mbox{-a.s.},
\end{equation}
then
%
\begin{equation}
\label{baselineactionmeasure} \prod_{V \in\X}
Z^V_0 \cdot\theta Q|_{\F_0}
\end{equation}
defines a counterfactual distribution on $\F_0$ that is absolutely
continuous with
respect to $P|_{\F_0}$ and imposes contemporaneously independent
outcomes.
\end{theorem}

\begin{pf}
First note that for every bounded $\F_0$-measurable random variable
$\eta$,
\begin{eqnarray*}
E_P \biggl[ \eta \frac{d \theta Q |_{\F_0}}{ d Q |_{\F_0}} \prod
_{V \in
\A} \frac{1 }{ Z_0^V } \biggr] & =& E_Q \biggl[
\eta \frac{d \theta
Q |_{\F_0}}{ d Q |_{\F_0}} \prod_{V \in
\X}
Z_0^V \biggr]
= E_{\theta Q} \biggl[ \eta \prod_{V \in
\X}
Z_0^V \biggr]
\\
& \leq& E_P [ \eta K ].
\end{eqnarray*}
This shows that \eqref{baselineactionmeasure} defines a finite
measure $P_\theta$ on $\F_0$ such that $P_\theta\ll P|_{\F_0}$.

We choose an enumeration $V_1, \ldots, V_m$ of the variables in
$\X$ such that $j < k$ implies that $T(V_j) \leq T(V_k)$.
If $V_k \in\X$ and $\eta$ is a bounded $\F_0^{
\{V_k
\} \cup p(V_k)}
$-measurable random variable, then
%
\begin{equation}
\label{blabla} E_Q \bigl[ \theta^* \eta | \F_0^{ \{ V_1, \ldots, V_{k-1} \}
}
\bigr] = \theta^* E_Q \bigl[ \eta | \F_0^{p(V)}
\bigr],\qquad  Q\mbox{-a.s.}
\end{equation}
To see this, we let $\eta_1$ be an $\F_0^{V_k}$-measurable
and bounded random variable and let $\eta_2$ be an $\F_0^{p(V_k)}$-measurable and
bounded variable and compute
\begin{eqnarray*}
E_Q \bigl[ \theta^* ( \eta_1 \eta_2) |
\F_0^{ \{ V_1,
\ldots, V_{k-1} \} } \bigr] & =& E_Q \bigl[
\eta_1 | \F_0^{ \{ V_1, \ldots, V_{k-1} \}} \bigr]\theta^*
\eta_2
\\
&=& \theta^* \bigl( E_Q \bigl[ \eta_1 |
\F_0^{ \{ V_1, \ldots,
V_{k-1} \} } \bigr]\eta_2 \bigr)
\\
&=& \theta^* \bigl( E_Q \bigl[ \eta_1 |
\F_0^{p(V_k) } \bigr]\eta_2 \bigr)
\\
&=& \theta^* E_Q \bigl[ \eta_1 \eta_2 |
\F_0^{p(V_k) } \bigr],\qquad  Q\mbox{-a.s.}
\end{eqnarray*}
Equation \eqref{blabla} now follows from the monotone class lemma.
Especially, this means that for every $k \leq m$,
%
\begin{equation}
E_Q \bigl[ \theta^* Z_0^{V_k} |
\F_0^{ \{ V_1, \ldots, V_{k-1}
\} } \bigr] = \theta^* E_Q \bigl[
Z_0^ {V_k} | \F_0^{p(V_k)} \bigr] = 1 ,\qquad Q
\mbox{-a.s.}
\end{equation}
and
\begin{eqnarray*}
E_{\theta Q} \bigl[ Z^{V_1}_0 \cdots
Z^{V_k}_0 \bigr] & =& E_{ Q} \bigl[ \theta^*
Z^{V_1}_0 \cdots\theta^* Z^{V_{k-1}}_0
E_Q \bigl[ \theta^* Z_0^{V_k} |
\F_0^{ \{ V_1, \ldots, V_{k-1} \}} \bigr] \bigr]
\\
& =& E_{ Q} \bigl[ \theta^* Z^{V_1}_0 \cdots
\theta^* Z^{V_{k-1}}_0 \bigr] = E_{ \theta Q} \bigl[
Z^{V_1}_0 \cdots Z^{V_{k-1}}_0 \bigr].
\end{eqnarray*}
That $P_\theta$
defines a probability measure on $\F_0$ follows by induction.

To see that \eqref{baselineaction1} and \eqref{baselineaction2} are
satisfied, suppose $V_k \in
\X$, and let $\eta, h$ be bounded
random variables such that $\eta$ is $\F_0^{V_k}$-measurable and
$h $ is $\F_0^{p(V_k)}$-measurable.
We see that
\begin{eqnarray*}
E_{P_\theta} [ \eta h ] & =& E_{ \theta Q} \Biggl[ \Biggl( \prod
_{j = 1}^{k-1} Z_0^{V_j}
\Biggr) \eta h Z_0^{V_j} \Biggr]
\\
& = &E_{ Q} \Biggl[ \Biggl( \prod_{j = 1}^{k-1}
\theta^* Z_0^{V_j} \Biggr) \theta^* h E_Q
\bigl[ \theta^* \eta Z_0^{V_k} | \F_0^{ \{ V_1, \ldots,
V_{k-1} \} }
\bigr] \Biggr]
\\
& =& E_{ Q} \Biggl[ \Biggl( \prod_{j = 1}^{k-1}
\theta^* Z_0^{V_j} \Biggr) \theta^* h \theta^*
E_Q \bigl[ \eta Z_0^{V_k} |
\F_0^{p(V_k)} \bigr] \Biggr]
\\
& =& E_{P_\theta} \bigl[ h \theta^* E_P \bigl[ \eta |
\F_0^{p(V)} \bigr] \bigr].
\end{eqnarray*}
If $V_k \in\A$, then
\begin{eqnarray*}
E_{P_\theta} [ \eta h ] & =& E_{ \theta Q} \Biggl[ \Biggl( \prod
_{j = 1}^{k-1} Z_0^{V_j}
\Biggr) \eta h Z_0^{V_j} \Biggr]
\\
& = &E_{ Q} \Biggl[ \Biggl( \prod_{j = 1}^{k-1}
\theta^* Z_0^{V_j} \Biggr) \theta^* h \theta^* \eta
E_Q \bigl[ Z_0^{V_k} | \F_0^{ \{ V_1, \ldots,
V_{k-1} \}}
\bigr] \Biggr]
\\
& =& E_{ Q} \Biggl[ \Biggl( \prod_{j = 1}^{k-1}
\theta^* Z_0^{V_j} \Biggr) \theta^* h \theta^* \eta \Biggr]
\\
& =& E_{P_\theta} \bigl[ h \theta^* \eta \bigr].
\end{eqnarray*}
\upqed\end{pf}
%

\subsection{Construction for the follow-up period}

Condition \eqref{action1} can be
made somewhat more concrete if
the processes, that may be directly intervened on, only
are allowed to jump at a given finite sequence of predictable times.
This behavior is very different from that of Poisson processes.
More formally, we assume that there exists
a bounded and $\F_t$-predictable multivariate counting
measure $\tilde U^A$ on $[0, T] \times J_A$
such that
%
\begin{equation}
\label{restrictedmodel} N |_{[0, T] \times J_A} \ll\tilde U^A
\end{equation}
for every $A \in\A$.
We can now show the reference measure $Q$ is quasi-invariant if the
probability of an outcome that coincides with the
counterfactually enforced outcome at short-term is not too small.
%
\begin{proposition} \label{restrictedquasiinvariance}
Suppose that $\theta$ is an $\F$-measurable transformation on
$\Omega$ that satisfies \eqref{action2}--\eqref{action5} and assume
\eqref{restrictedmodel}.
If there exists a bounded and $\PP$-measurable process $\tilde Y$
such that:
\begin{longlist}[(1)]
\item[(1)]
$\theta Q |_{\F_0} \ll Q|_{\F_0} $;
\item[(2)]
%
\begin{equation}
\int_{J_A}\int_0^T h(s,
x) \theta^* N (ds, dx) =\int_{J_A}\int_0^T
h(s, x) \tilde Y(s, x) U^A (ds, dx)
\end{equation}
$Q$-a.s. for every $A \in\A$ and bounded and $\PP$-measurable process $h$;
\item[(3)]
there exists a constant $c > 0$ such that
%
\begin{equation}
1 - \theta^* N\bigl( \{s\}, J_A \bigr) \leq c\cdot \bigl(1 -
U^A\bigl( \{s\}, J_A \bigr) \bigr),\qquad Q\mbox{-a.s. }
\end{equation}
for every $s \in[0, T]$,
\end{longlist}
then $\theta Q \ll Q$.\vadjust{\goodbreak}
\end{proposition}

\begin{pf}

The integral equation
\begin{eqnarray*}
& &\int_J \int_0^T h(s,
x) U^\theta(ds, dx)
\\
&&\qquad=  \sum_{A \in\A} \int_{J_A} \int
_0^T h(s, x) \theta^* N (ds, dx) + \sum
_{V \in\X} \int_{J_A} \int
_0^T h(s, x) U^V (ds, dx)
\end{eqnarray*}
defines an $\F_t$- predictable random measure $U^\theta$ on $[0, T]
\times J$.

Let $B \subset J$ be a measurable subset, and define
$N_t^B := \int_0^t \int_{B}N(ds, dx)$. If
$B \subset J_A$ for an $A \in\A$ and $S$ is
a $\F_t$-adapted stopping time, then
\[
E_{\theta Q} \bigl[ N_S^B - U^\theta_S
\bigl( B , [0, t] \bigr) \bigr] = E_{\theta Q} \bigl[ N_S^B
- \theta^*N_S^B \bigr] = E_{\theta Q} \bigl[
\theta^* N_S^B - \theta^*N_S^B
\bigr] = 0.
\]
This means that $N_t - U^\theta_t( B , [0, t] )$ defines a local
$Q$-martingale with respect to the filtration $\F_t$. Similarly, if $B
\subset J_\X$, note that
\begin{eqnarray*}
E_{\theta Q} \biggl[ \int_0^T
h_s \,dN_s^B \biggr] &= & E_{ Q}
\biggl[ \int_0^T \theta^* h_s
\,dN_s^B \biggr]
\\
&= & E_{ Q} \biggl[ \int_0^T
\theta^* h_s \,dU( ds, B) \biggr]
\\
&= & E_{ \theta Q} \biggl[ \int_0^T
h_s \,dU^\theta( ds, B) \biggr]
\end{eqnarray*}
for every bounded and $\F$-predictable process $h$.
Now, $N( [0, t], B ) - U^\theta( [0, t],\break B)$ defines a local
$\theta Q$-martingale with respect to the filtration $\{\F_t\}_{t
\in[0,T]} $.
This means that
\[
E_{\theta Q} \biggl[ \int_J \int
_0^T h(s,x) N(ds, dx) \biggr] =
E_{\theta Q} \biggl[ \int_J \int_0^T
h(s,x) U^\theta(ds, dx) \biggr]
\]
for every bounded and $\PP\otimes\J$-measurable process $h$.

We define the processes
\begin{eqnarray*}
H^A(t, x) &:=& \tilde Y(t, x) - 1 - \frac{U( \{t\}
, J_A ) - \theta^* N( \{t\} , J_A) }{ 1- U( \{t\}
, J_A ) } I\bigl( U\bigl(\{t\}, J_A\bigr) \neq1\bigr),
\\
\zeta^A_t & :=& \int_{J_A} \int
_0^t H^A(s, x) \bigl( N( ds, dx) -
U (ds, dx) \bigr),
\end{eqnarray*}
and let $\zeta:= \sum_{A \in\A} \zeta^A$.

By \cite{JacodShiryaev}, Proposition~I 3.13, there
exists a $\PP$-measurable and nonnegative stochastic process
$\gamma^A $ such that $\gamma^A \leq1$ and
\[
\int_{J_A} \int_0^T h(s,
x) U^A (ds, dx) = \int_{J_A} \int
_0^T h(s, x) \gamma^A (s,x) \tilde
U^A (ds, dx)
\]
$Q$-a.s.
for every bounded and $\PP$-measurable stochastic process $h$.\vadjust{\goodbreak}

A computation shows that the predictable variation process for $\zeta$ with
respect to $Q$ satisfies
\begin{eqnarray*}
\langle\zeta, \zeta\rangle_t &= &\sum_{A \in\A}
\bigl\langle\zeta^A, \zeta^A \bigr\rangle_t
\\
&= & \sum_{A \in\A} \int_{J_A} \int
_0^t H^A(s, x)^2
\gamma^A(s, x) \bigl( 1- \gamma^A(s, x)\bigr) \tilde
U^A( ds, dx),
\end{eqnarray*}
which is $Q$-a.s. uniformly bounded.
Now,
\cite{LepingleMemin}, Theorem~II.1, implies that the SDE
%
\begin{equation}
\rho_t = \frac{d\theta Q|_{\F_0} } {
d Q|_{\F_0} } + \int_0^t
\rho_{s-} \,d \zeta_s
\end{equation}
defines a uniformly integrable $Q$-martingale with respect to the
filtration $\F_t$. This means that
\[
\tilde Q:= \rho_T \cdot Q
\]
defines a
probability measure on $( \Omega, \F)$.

A computation shows that if $B \subset J_V$ for some $V \in\V$, then
%
\begin{equation}\quad
N_t^B - U_t\bigl( [0, t], B \bigr) - \int
_0^t \rho^{-1}_{s-} \,d
\bigl\langle N^B - U^B , \rho \bigr\rangle_s
= N_t^B - U^\theta\bigl( [0, t] , B\bigr).
\end{equation}
Girsanov's Theorem \cite{JacodShiryaev}, Theorem~III 1.21, implies that
\[
E_{\tilde Q} \biggl[ \int_J \int
_0^T h(s,x) N(ds, dx) \biggr] =
E_{\tilde Q} \biggl[ \int_J \int_0^T
h(s,x) U^\theta(ds, dx) \biggr]
\]
for every bounded and $\PP\otimes\J$-measurable process $h$.
Finally, \cite{JacodMultivariate}, Theorem~3.4, implies that there
exists only one probability measure which has $U^\theta$ as a dual
predictable projection for $N$. Therefore $\theta Q = \tilde Q \ll
Q$.
\end{pf}

The next result is important and says that
if the probability of observing an outcome that coincides with the
counterfactually enforced outcome at short-term is not too small,
then there exists a counterfactual distribution for the follow-up
period.
The counterfactual distribution can then be obtained by re-weighting
the
factual distribution, that is, $P_\theta\ll P$.
Note that \eqref{finfaktorisering} provides a continuous-time
analogy of the truncated factorization formula for Bayesian
networks \cite{Pearl}, (3.10).

\begin{theorem} \label{hovedeksistens}
Suppose that the conditions of Theorem~\ref{baselinexistthm}
are satisfied and that
there exists a bounded and $\PP$-measurable process $Y$
such that:
\begin{longlist}[(1)]
\item[(1)]
%
\begin{equation}
\qquad\int_{J_A}\int_0^T h(s,
x) \theta^* N (ds, dx) =\int_{J_A}\int_0^T
h(s, x) Y(s, x) \lambda(s,x)U (ds, dx)
\end{equation}
$P$-a.s. for every $A \in\A$ and bounded and $\PP$-measurable process $h$;\vadjust{\goodbreak}
\item[(2)]
there exists a constant $c > 0$ such that
%
\begin{equation}
1 - \theta^* N\bigl( \{s\}, J_A \bigr) \leq c \bigl(1 - \lambda
\cdot U\bigl( \{s\}, J_A \bigr) \bigr) ,\qquad P\mbox{-a.s. }
\end{equation}
for every $s \in[0, T]$.
\end{longlist}

Then there exists a counterfactual distribution $P_\theta$ such that
$P_\theta\ll
P$.
We also have
that $P_\theta\ll\theta Q$ and
%
\begin{equation}
\label{finfaktorisering} X_t := \prod_{V \in\X}
Z^V_t,
\end{equation}
where $Z^V$ is the process defined in \eqref{ZV},
defines a $\theta Q$-martingale with respect to the filtration
$\{\F_t\}$ that
satisfies the SDE
%
\begin{equation}
\label{Xequation} X_t = \prod_{V \in\mathcal X}
Z_0^V + \sum_{V \in\mathcal X} \int
_0^t X_{s-} \,d K^V_s
\end{equation}
and
%
\begin{equation}
\frac{d P_\theta }{d \theta Q} = X_T.
\end{equation}
\end{theorem}

\begin{pf}
We follow the proof of Proposition~\ref{restrictedquasiinvariance} and
define the processes
\begin{eqnarray*}
G^A(t, x) &:=& Y(t, x) - 1 - \frac{\lambda\cdot U( \{t\}
, J_A ) - \theta^* N( \{t\} , J_A) }{ 1- \lambda\cdot U( \{t\}
, J_A ) } I\bigl( \lambda\cdot
U\bigl(\{t\}, J_A\bigr) \neq1\bigr),
\\
\xi^A_t & :=& \int_{J_A} \int
_0^t G^A(s, x) \bigl( N( ds, dx) -
\lambda\cdot U (ds, dx) \bigr),
\\
\xi& :=& \sum_{A \in\A} \xi^A.
\end{eqnarray*}

By \cite{JacodShiryaev}, Proposition~I 3.13, there
exists a $\PP$-measurable and nonnegative stochastic process
$\gamma^A $ such that $\gamma^A \leq1$ and
\[
\int_{J_A} \int_0^T h(s,
x) \lambda(s,x) U (ds, dx) = \int_{J_A} \int
_0^T h(s, x) \gamma^A (s,x) \tilde
U^A (ds, dx)
\]
$Q$-a.s.
for every bounded and $\PP$-measurable stochastic process $h$.

A computation shows that the predictable variation process for $\xi$ with
respect to $P$ satisfies
\begin{eqnarray*}
 \langle\xi, \xi\rangle_t &= &\sum_{A \in\A}
\bigl\langle\xi^A, \xi^A \bigr\rangle_t
\\
&= & \sum_{A \in\A} \int_{J_A} \int
_0^t G^A(s, x)^2
\gamma^A(s, x) \bigl( 1- \gamma^A(s, x)\bigr) \tilde
U^A( ds, dx),
\end{eqnarray*}
which is $Q$-a.s uniformly bounded.
Now,
\cite{LepingleMemin}, Theorem~II.1, implies that the SDE
%
\begin{equation}
\label{Wdef} W_t = \frac{dP_\theta|_{\F_0} } {
d P|_{\F_0} } + \int_0^t
W_{s-} \,d \xi_s
\end{equation}
defines a uniformly integrable $P$-martingale with respect to the
filtration $\F_t$. This means that
\[
P_\theta:= Z_T \cdot P
\]
defines a
probability measure on $( \Omega, \F)$.

The integral equation
%
\begin{eqnarray}
\label{nutheta}\qquad && \int_J \int_0^T
h(s , x) \nu^\theta( ds, dx)
\\
&&\qquad =  \int_{J_\X} \int_0^T
h(s , x) \lambda(s,x) U( ds, dx) + \int_{J_\A} \int
_0^T h(s , x) \theta^* N ( ds, dx)
\end{eqnarray}
defines a predictable and nonnegative random measure $\nu^\theta$ on
$[0, T] \times J$
such that
\begin{eqnarray*}
\xi_t &=& \int_J \int_0^t
 \biggl( Y(s, x) - 1 - \frac{U \lambda\cdot(\{s\}, J ) - \nu^\theta(\{s\}, J ) } {
1- \lambda\cdot U(\{s\}, J )} \\
&&\hspace*{111pt}{}\times I \bigl( \lambda\cdot U \bigl(\{s\}, J
\bigr) \neq1 \bigr) \biggr)
 N(ds, dx) - \lambda\cdot U ( ds, dx).
\end{eqnarray*}
We obtain from \cite{JacodMultivariate}, Theorem~5.2, that
\[
E_{P_\theta} \biggl[ \int_J \int
_0^T h(s,x) N(ds, dx) \biggr] =
E_{P_\theta} \biggl[ \int_J \int_0^T
h(s,x) \nu^\theta(ds, dx) \biggr]
\]
for every bounded and $\PP\otimes\J$-measurable process $h$; that is,
$P_\theta$ defines a counterfactual distribution.

We may compute that
\begin{eqnarray*}
\Delta\zeta_s^A &= & \int_{J_A}
\tilde Y(s, x) N\bigl( \{s\} , dx\bigr) - \theta^* N\bigl( \{s\}, J_A
\bigr)
\\
&&{} + \bigl( \tilde U \bigl( \{s\}, J_A\bigr)- \theta^* N\bigl( \{s
\}, J_A\bigr) I \bigl( U \bigl( \{s\} , J_A\bigr) \neq1
\bigr) \bigr)
\\
&&\quad{} \times \bigl( \tilde U \bigl( \{s\}, J_A\bigr) - N\bigl( \{s\},
J_A\bigr) \bigr)
\end{eqnarray*}
and that
\begin{eqnarray*}
\Delta\xi_s^A &= & \int_{J_A} Y(s,
x) N\bigl( \{s\} , dx\bigr) - \theta^* N\bigl( \{s\}, J_A\bigr)
\\
& &{}+ \bigl( \tilde U \bigl( \{s\}, J_A\bigr)- \theta^* N\bigl( \{s
\}, J_A\bigr) I \bigl( \lambda\cdot U \bigl( \{s\} , J_A
\bigr) \neq1 \bigr) \bigr)
\\
&&\quad{} \times \bigl( \tilde U \bigl( \{s\}, J_A\bigr) - N\bigl( \{s\},
J_A\bigr) \bigr).
\end{eqnarray*}
We moreover define a process $\chi$ as follows:
\[
\chi_t := \sum_{s \leq t} \frac{\Delta\xi_s - \Delta\zeta_s } {
\Delta\zeta_s + 1 }
I ( \Delta\zeta_s \neq-1 ).
\]
One can show that $\chi$ only jumps at the jump times of $\tilde U$
and that $\Delta\chi$ is uniformly bounded.
This means that the SDE
%
\begin{equation}
\pi_t := \frac{d\theta Q |_{\F_0} }{ d Q |_{\F_0} } \prod_{A
\in\mathcal A}
\frac{1}{Z_0^A} + \int_0^t
\pi_{s-} \,d \chi_s
\end{equation}
defines a $P$ semi-martingale with respect to the filtration $\F_t$.
Note that $\Delta\zeta_s = -1$ implies that
$\Delta\xi_s = -1$, so
%
\begin{equation}
\zeta+ [ \zeta, \chi] + \chi= \xi, \qquad P\mbox{-a.s.}
\end{equation}
Yor's additive formula \cite{Protter}, Theorem~II 38, then implies
that
%
\begin{equation}
\pi_t \rho_t = \frac{d P_\theta|_{\F_0}} {d P |_{\F_0}} + \int
_0^t \pi_{s-} \rho_{s-} \,d
\xi_s.
\end{equation}
This implies that $W = \rho\pi$, and hence
\[
E_{P_\theta} [h ] = E_{P} [h W_T ] =
E_{Q} [h Z_T \rho_T \pi_T ] =
E_{\theta Q} [h Z_T \pi_T ]
\]
for every bounded and $\F_T$-measurable random variable $h$,
so $P_\theta\ll \theta Q$.
Finally~\cite{JacodMultivariate}, Theorem~5.1, shows that the
likelihood ratio\vspace*{1pt} $\frac{d
P_\theta}{d\theta Q}$ is given by the SDE~\eqref{Xequation}, and
hence Yor's additive formula provides identity
\eqref{finfaktorisering}.
\end{pf}

Note that since $P_\theta\ll\theta Q = \theta^2 Q$, the
counterfactual distribution $P_\theta$
is actually invariant with
respect to the action $\theta$, that is,
%
\begin{equation}
\label{actioninvariance} \theta P_\theta= P_\theta.
\end{equation}

\section{Local independence} \label{LI}
\subsection{Identifiability and short-term dependence}
A causal effect is identifiable if it can be uniquely obtained
from the factual distribution of the observable variables. This is
generally very hard to determine and may also require further
parametric assumptions.
We show
that it is possible to take advantage of graphical structure, in
terms of local independence graphs, to do this.
Such graphs are useful when deciding in which
situations causal effects are identifiable, and
also which factors we might adjust for.

We will say that $V \in\V$ is \textit{locally independent} of
a subset $\mathcal B \subset\V$
at baseline,
conditionally on $\V' \subset\V$, if the conditional density of
$V_0$, given the past, does not depend on the baseline
information from $\mathcal B$. More precisely,
for every integrable and $\F_0^V$-measurable random variable $\eta$,
there exists a random variable $\tilde\eta$ that is\vadjust{\goodbreak}
$ \F_0^{ p(V) \cap( \V' \setminus
\mathcal B) }$-measurable and such that if $h$ is
$\F_0^{p(V) \cap \V' }$-measurable, then
%
\begin{equation}
E_P [\eta h ] = E_P [\tilde\eta h ].
\end{equation}

A process $V \in\V$ is locally independent of $\mathcal B \subset
\V$ during follow-up,
conditionally on $\V'$,
if
for every process $X$ on the form
\eqref{detectableDef},
there exists an
$\F^{\{V\} \cup\V' \setminus\mathcal B}_t$-predictable process
$\Lambda$ with finite variation
such that
%
\begin{equation}
\label{condindfollow-up} E_P \biggl[ \int_0^Th_s
\,dX_s \biggr] = E_P \biggl[ \int_0^Th_s
\,d\Lambda_s \biggr]
\end{equation}
for every bounded and $\F^{\{V\} \cup\V' }_t$-predictable process
$h$.
If
$V$ is locally independent of $\mathcal B$,
conditionally on $\V'$, both at baseline and during follow-up, we will
say that $V$ is locally independent of $\mathcal B$,
conditionally on $\V'$.
This will sometimes be written $\mathcal B
\nrightarrow V | \V'$.
A \textit{local independence graph} is a directed graph
$G = (\V', \mathcal E)$ for $\V' \subset\V$ such that the absence
of an arrow
from a subset $\mathcal B \subset\V'$ to a process $V \in\V'$
means that
$\mathcal B \nrightarrow V | \V'$. Note that local independence graphs
are also
refered to as local independence graphs (see \cite{DidelezRoyalSB,AalenRoyslandRSSA}) and were introduced in \cite{Schweder}.

Given time points $\{T(V_0)\}_{V \in\V}$ at baseline
and a local independence graph $G = (\V, \E)$, we can pick a linear
ordering of $\V_0$ that satisfies
\eqref{ordering} and therefore yields
%
\begin{equation}
\label{ordereddirectedMarkovProperty} V^i_0
\ci_P \bigl\{V_0^1, V_0^2
, \ldots, V^{i-1}_0 \bigr\} | \F_0^{\mathrm{pa}(V^i)}
\end{equation}
for every $i \leq n$. Property \eqref{ordereddirectedMarkovProperty} is known as the ordered directed
Markov Property and was shown to be equivalent to the local
directed Markov property in \cite{Lauritzen}, Theorem~2.11. This
means that Bayesian networks and local independence graphs are two
descriptions of the same structure when the nodes correspond to
single variables. Note that local independence graphs, where the
nodes are allowed to be families of variables or processes, are
allowed to be cyclic.

\subsection{Measurability of intensities}

Local independence during the follow-up is closely related to the
measurability of intensities.
\def\B{\mathcal B}
%
\begin{lemma} \label{elocind}
Suppose that $V$ is locally independent of $\mathcal B $
at baseline,
conditionally on $\V'$, then
$ \B\nrightarrow V | \V'$ if and only if there exists
a nonnegative and $\PP^{\{V\} \cup\V' \setminus\B}$-measurable
process $\lambda^V$
such that
%
\begin{equation}
\label{explicitlocind} E_P \biggl[ \int_{J_V}
\int_0^T h(s,x)N(ds, dx) \biggr] =
E_P \biggl[ \int_{J_V} \int_0^T
h(s,x) \lambda^V(s, x) U(ds, dx) \biggr]\hspace*{-35pt}
\end{equation}
for every bounded and $\PP^{\{V\}\cup\V'}$-measurable process
$h$.
\end{lemma}

\begin{pf}
If there exists a process $\lambda^V$ as in \eqref{explicitlocind},
then $\B\nrightarrow V |\V'$ follows directly. Conversely,
suppose that $\B\nrightarrow V |\V'$ and let $D \subset J_V$
be a measurable subset. Now,\vadjust{\goodbreak}
$
N^D_t := N( [0, t] ,D )
$
defines a processes on the form \eqref{detectableDef}, so there must
exist a corresponding
predictable increasing process $\Lambda^D$ of finite variation such that
\[
E_P \biggl[ \int_0^T
h_s \,dN_s^D \biggr] = E_P
\biggl[ \int_0^T h_s d
\Lambda_s^D \biggr]
\]
for every bounded and $\F_t^{ \{V\} \cup\V' }$-predictable
process $h$.

The Radon--Nikodym theorem now provides an $\F^{ \{V\} \cup\V'
\setminus\B }$-measurable and nonnegative process $\lambda^{(D)}$
such that
%
\begin{equation}
E_P \biggl[ \int_0^T
h_s \,dN_s^D \biggr] = E_P
\biggl[ \int_0^T h_s
\lambda_s^{(D)} U( ds, D ) \biggr]
\end{equation}
for every bounded and $\F_t^{\{V\} \cup\V' }$-measurable process $h$.

Since $J$ is a Lusin space,
we may construct a nonnegative and
$\PP^{\{V\} \cup\V' }$-measurable process $\lambda^V$ that satisfies
\eqref{explicitlocind}
as a limit of
processes that are finite linear combinations of processes on the form
$f \cdot J_D$, where $D$ is a measurable subset in $J_V$, and $f$ is
a bounded $\F_t^{\{V\} \cup\V' \setminus\B }$-measurable process.
\end{pf}
%

\subsection{Markovian factorization property}

The local Markov property implies the \textit{Markovian factorization
property}; see
\cite{Pearl}, (1.33) and \cite{Lauritzen}, (2.10).
We
will now see that a local independence graph yields a similar
factorization for the follow-up period. We use the following notation from
graph theory: whenever $V \in\V$, let $cl(V) \subset\V$
denote the set formed by $V$ and its parents in $G$.
%
\begin{theorem} \label{factorizationoflikelihood}
If $G = ( \V, \mathcal E)$ is a local independence graph with respect
to $P$, then
there
exists an $\F_t^{cl(V)}$-adapted $P$-indistinguishable
version of each process $Z^V$ from
Theorem~\ref{likfactorization} where
\[
Z = \prod_{V \in\V} Z^V, \qquad P\mbox{-a.s.}
\]
\end{theorem}

\begin{pf}
Let $\F_0^{\mathrm{pa}(V)} := \bigvee_{V' \in\operatorname{pa}(V)} \F_0^{V'}$ and
$\F_0^{\mathrm{cl}(V)}:= \bigvee_{V' \in\operatorname{cl}(V)} \F_0^{V'}$ and
let
\[
Y^V := \frac{ dP|_{\F_0^{\mathrm{pa}(V)} } } {
dQ|_{\F_0^{\mathrm{pa}(V)} } }.
\]
Now
%
\begin{equation}
P|_{\F_0^{\mathrm{cl}(V)} } \ll Y^V \cdot Q|_{\F_0^{\mathrm{cl}(V)} },
\end{equation}
so there exists, by the Radon--Nikodym theorem, an
$\F_0^{\mathrm{cl}(V)}$-measurable random variable $\tilde Z_0^V$ such
that
%
\begin{equation}
P|_{\F_0^{\mathrm{cl}(V)}} = \tilde Z_0^V Y^V \cdot
Q|_{\F_0^{\mathrm{cl}(V)}}.\vadjust{\goodbreak}
\end{equation}
We then have, for every bounded and measurable function $h$, that
\begin{eqnarray*}
E_P \bigl[ h(V_0) | \F_0^{p(V)}
\bigr] & =& E_P \bigl[ h(V_0) | \F_0^{pa(V)}
\bigr] = E_Q \bigl[ h(V_0) \tilde Z_0^V
| \F_0^{pa(V)} \bigr]
\\
& =& E_Q \bigl[ h(V_0) \tilde Z_0^V
| \F_0^{p(V)} \bigr].
\end{eqnarray*}
The contemporaneous independence at baseline and a simple monotone class
argument shows that
%
\begin{equation}
E_P [ \eta ] = E_Z \biggl[ \eta\prod
_{V \in\V} \tilde Z_0^V \biggr]
\end{equation}
for every bounded and $\F_0$-measurable random variable $\eta$.

For the follow-up, note that by Lemma~\ref{elocind} there exists
a nonnegative and $\PP^{\mathrm{cl}(V) }$-measurable process $\lambda^V$
such that
\[
E_P \biggl[ \int_{J_V} \int
_0^T h(s,x)N(ds, dx) \biggr] = E_P
\biggl[ \int_{J_V} \int_0^T
h(s,x) \lambda^V(s, x) U(ds, dx) \biggr]
\]
for every bounded and $\PP$-measurable process $h$.

We may now form $K^V$, $Z^V$
and $Z$ as in Theorem~\ref{likfactorization} using $\lambda^V$
instead of $\lambda$.
Following the short argument in \cite{Bremaud}, Theorem~II T12, we see
that any other choice of a nonnegative and
$\PP$-measurable process $\lambda$ that satisfies the
previous equation would necessarily give
%
\begin{equation}
\int_{J_V} \int_0^T I
\bigl( \lambda(s,x) \neq\lambda^V(s,x) \bigr) N(ds, dx) = 0,\qquad P\mbox{-a.s.}
\end{equation}
This means that the corresponding versions of the process $K^V$
from \eqref{KV}
would be~$P$ indistinguishable. Furthermore, this also means that the version
corresponding to $\lambda^V$ provides an $\F_t^{\mathrm{cl}(V)}$-adapted
solution of the SDE \eqref{ZV} which is
$P$-indistinguishable version from $Z^V$.
\end{pf}
%

\section{An example: Controlled direct effects} \label{CE}
We now illustrate how local independence graphs can be used to
identify causal effects by an example with cancer patients.
Suppose each
patient is offered one of two different surgical treatments, $a_1$ or $a_2$.
The patient is subject to an examination after surgery where some
measurements are taken.
These measurements might depend on the chosen
surgical procedure and some underlying health condition that is not
directly observed.
After the surgery, the patient is given further treatment in order to
prevent relapse. The chosen post surgery treatment strategy might
depend on the
surgical procedure and the measurements.

We consider a generic model for the patients in this scenario. The
relevant outcomes are provided by the family of random variables
$\V= \{W, A , L, K , B\}$. As in Section~\ref{observationalregime},
we consider a probability measure $Q$ such that these variables are
independent and a probability measure $P$ that governs the frequency
of outcomes in the factual scenario and such that $P \ll Q$.
Let the random variable $A$ denote the choice of surgery, let $W$ denote
the latent health condition, let $L$ take the value of the
measurements after surgery, let $K$ denote the post surgery
treatment strategy and let $B$ denote the status of relapse.
We furthermore assume that $T(W) < T(A) < T(L) < T( K) <T( B)$ and
that the following local independencies are satisfied:\vspace*{4pt}
\begin{center}
\includegraphics{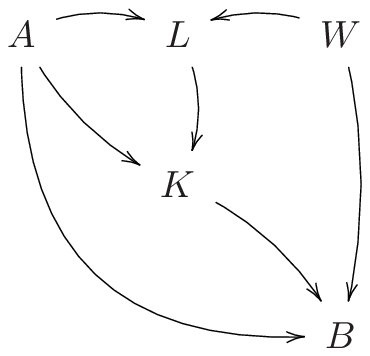}\vspace*{4pt}
\end{center}

How much of the treatment
effect is due
to the choice of surgical procedure alone, that is, not due to the choice
of
post surgery treatment?
Pearl
\cite{Pearl}, Section~4.5.3, showed that it is possible to identify the
controlled direct effect from surgery on the risk of relapse, even
without any observations of $W$. We rephrase his argument slightly:
%
\begin{proposition}
If $\theta^* K$ is $\F_0^L$-measurable, $\theta^* A$ is constant,
$L$, $W$ and $B$ are $\theta$-invariant, there exists a constant $c >
0$
such that
%
\begin{equation}
\label{concreteassumption} P \bigl( A = \theta^* A \bigr) > 0\quad \mbox{and}\quad P
\bigl( K = \theta^* K | A = \theta^* A , L \bigr) \geq c,\qquad P\mbox{-a.s}
\end{equation}
and $h$ is a bounded and measurable function, then there exists a
unique counterfactual distribution $P_\theta$ such that $P_\theta\ll
P$ and
%
\begin{equation}
E_{P_\theta} \bigl[ h(B) \bigr] = \theta^* E_P \bigl[ \theta^*
E_P \bigl[ h(B) | \F_0^{ \{ L, A, K \}} \bigr] |
\F_0^A \bigr], \qquad P_\theta\mbox{-a.s. }
\end{equation}

Let $\tilde \F_0 := \F_0^{ \{ L, A, K , B \} } $ and suppose that
$\tilde Z^B$ is a nonnegative and $\tilde
\F_0$-measurable random
variable and $\tilde Z^L$ is a nonnegative $\F_0^{A}$-measurable
random variable such that
\begin{eqnarray*}
E_P \bigl[ h(B) | \F_0^{ \{A, L , K\}} \bigr] & =&
E_Q \bigl[ h(B) \tilde Z^B | \F_0^{ \{A, L , K\}}
\bigr],
\\
E_P \bigl[ h(L) | \F_0^{A} \bigr] & =&
E_Q \bigl[ h(L) \tilde Z^L | \F_0^{A}
\bigr]
\end{eqnarray*}
$P$-a.s.
Now,
%
\begin{equation}
\label{Gidentify} E_{P_\theta} [ H ] = E_{\theta Q} \bigl[ H \tilde
Z^L \tilde Z^B \bigr]
\end{equation}
for every $\tilde\F_0$-measurable random variable $H$, that is,
%
\begin{equation}
\label{identifylikelihood} \frac{d P_\theta|_{ \tilde\F_0 }}{d \theta Q |_{ \tilde\F_0} } = \tilde Z^L \tilde
Z^B.
\end{equation}
\end{proposition}

\begin{pf}
Note that \eqref{concreteassumption} means
that \eqref{baselinecondition} is satisfied, that is, we obtain a
counterfactual distribution $P_\theta$ from
Theorem~\ref{baselinexistthm}.\vadjust{\goodbreak}

Whenever $h_1, h_2$ are bounded and measurable functions, then
\begin{eqnarray*}
E_{P_\theta} \bigl[h_1(B) h_2(L) \bigr] & =&
E_P \bigl[W_0 h_1( B) h_2( L)
\bigr]
\\
& =& E_P \bigl[ W_0 E_P \bigl[h(B) |
\F_0^{A, K , L} \bigr] h(L) \bigr]
\\
& =& E_{P_\theta} \bigl[ E_P \bigl[h(B) |\F_0^{A, K , L}
\bigr] h(L) \bigr]
\\
& =& E_{P_\theta} \bigl[ \theta^* E_P \bigl[h(B) |
\F_0^{A, K , L} \bigr] h(L) \bigr]\qquad \mbox{by
\eqref{actioninvariance}}.
\end{eqnarray*}
This shows that
$ E_{P_\theta}  [ h_1(B)  | \F_0^L  ] = \theta^* E_P
[
h_1(B)  | \F_0^ { \{A, L , K\}}  ]
$
$P_\theta$-a.s.
Moreover, note that
\begin{eqnarray*}
 E_{P_\theta} \bigl[ h_2 (L) \bigr]& =& E_{P_\theta}
\bigl[ \theta^* E_P \bigl[ h_2(L) |
\F_0^{A, W} \bigr] \bigr]
\\
&= & E_{P_\theta} \bigl[ \theta^* E_P \bigl[
h_2(L) | \F_0^{A} \bigr] \bigr] \\
&=& \theta^*
E_P \bigl[ h_2(L) | \F_0^{A}
\bigr],\qquad P_\theta\mbox{-a.s.}
\end{eqnarray*}

Combining these computations, we obtain
\begin{eqnarray*}
E_{P_\theta} \bigl[ h(B) \bigr] & =& E_{P_\theta} \bigl[
E_{P_\theta} \bigl[ h(B) | \F_0^L \bigr] \bigr]\\
& =&
E_{P_\theta} \bigl[ \theta^* E_P \bigl[ h(B) |
\F_0^{L , A , K} \bigr] \bigr]
\\
& =& \theta^* E_P \bigl[ \theta^* E_P \bigl[ h(B) |
\F_0^{L , A , K} \bigr] | A \bigr]
\end{eqnarray*}
$P_\theta$-a.s. for every bounded and measurable function $h$.

To see that equation
\eqref{Gidentify} is satisfied, note that by the monotone class
lemma,
\begin{eqnarray*}
E_{P_\theta} [ H ] & = &\theta^* E_P \bigl[ \theta^*
E_P \bigl[ H | \F_0^ { \{A, L , K\}} \bigr] |
\F_0^A \bigr]
\\
& =& \theta^* E_Q \bigl[ \theta^* E_Q \bigl[ H \tilde
Z^B | \F_0^ {
\{A, L , K\}} \bigr] \tilde
Z^L | \F_0^A \bigr]
\\
& =& E_{\theta Q} \bigl[ E_Q \bigl[ H \tilde Z^B |
\F_0^ { \{A, L
, K\}} \bigr] \tilde Z^L \bigr]
 \\
 &=& E_{\theta Q} \bigl[ H \tilde Z^B \tilde Z^ L
\bigr] .
\end{eqnarray*}
\upqed\end{pf}

If we consider actions $\theta_1$ and $\theta_2$ such that $\theta_i^* A = a_i$ and $\theta_1^* K = \theta_2^*
K$, $Q$-a.s. then the relative direct risk of relapse is given by
%
\begin{equation}
\frac{P_{\theta_1}  ( B = 1  ) } {P_{\theta_2}  ( B = 1
) } = \frac{ E_P  [ \theta^*_1 E_P
[ h(B) | \F_0^{ \{A, L , K\}}  ]  |A = a_1  ] } {
E_P  [ \theta^*_2 E_P
[ h(B) | \F_0^{ \{A, L , K\}}  ]  | A = a_2  ] } .
\end{equation}

\subsection{Incomplete observations and time dependent treatments}

We have not yet taken into account that the patient observations could
be censored during the follow-up period. There might be several
reasons for such censoring. This might be due to the end of study
period, drop-out due
to the underlying health or because of other reasons.
The risk of having an observed relapse will typically be smaller than
the risk of having a relapse.
We will work in the framework of event history analysis
in order to provide a reasonable effect
measure subject to such incomplete observations. This will also
allow us to consider time dependent post surgery strategies $K$.

\subsubsection{A dynamic model}
We proceed with the previous setup, but where
$B$ and $K$ are represented by processes and every patient
may be censored during the follow-up period.
The factors $A, L$ and $W$ are as in the previous example.
$B$ is represented by a counting process that jumps from $0$ to $1$ at
the time of the event. The censoring of the
individual is represented by a counting process $C$ that
jumps from $0$ to $1$ at the time of censoring.

We suppose that the baseline treatment $A$ may be of two different
types; hence $A$ takes value in $\{0, 1\}$.
Moreover, we suppose that additional post-surgery treatment is given
to the patient
at the jumps of the counting process $K$.
This treatment may be given recursively, but only
at a series of $\F_t$-predictable times; that is, \eqref{restrictedmodel} must be satisfied.
We furthermore suppose that $\theta^* K_s$ is constant for every $s$
$P$-a.s. and suppose that
$B_0 = 0$, $K_0 = 0$ and $C_0 = 0$ $P$-a.s.

Let $T_1, \ldots, T_n$ denote the potential post-treatment times, and
let $U_t^K := \sum_{i} I ( T_i \leq t)$. The counting process $U^K$
is predictable and
$\nu^K_t = \int_0^t P  ( \Delta K_s \neq0  | \F_{s-}  )
\,dU^K_s$. By Theorem~\ref{hovedeksistens}, we see that there exists a
counterfactual distribution if
$P  ( A = \theta^* A  ) > 0$
and there exist $c_1 , c_2 > 0$ such that
%
\begin{equation}
1- c_1 P ( \Delta K_s = 0 | \F_{s-} ) \leq
\Delta \theta^* K_s \leq c_2 P ( \Delta
K_s \neq0 | \F_{s-} )
\end{equation}
for every $s$ $P$-a.s.

We suppose that the following local independence graph is satisfied with
respect to the factual distribution $P$:\vspace*{4pt}
\begin{center}
\includegraphics{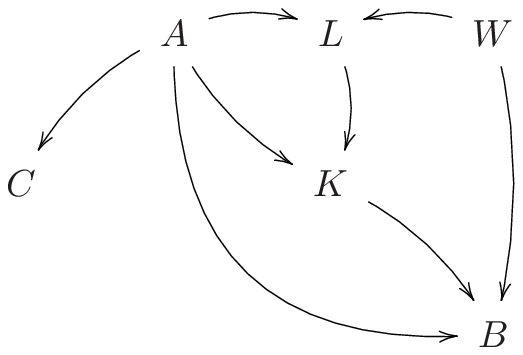}\vspace*{4pt}
\end{center}
Especially, this means that the short-term behavior of the censoring
may not depend on other variables than $A$.

\subsubsection{Restriction to Aalen's additive hazard model}
If we assume that the event process satisfies Aalens's additive
hazard model \cite{Andersen}, it is actually possible to identify, and also
consistently estimate\vadjust{\goodbreak}
the direct effect from surgery.
Every outcome after the time of censoring is supposed to be
unobserved. In addition, we assume that we are not able to
observe the variable~$W$.

We
consider the censored process
\begin{eqnarray*}
\tilde B_t := B_0 + \int_0^t
( 1- C_{s-}) \,dB_s
\end{eqnarray*}
and let
$\tilde\F_t$ denote the filtration that is generated by $A, K,
L, C $ and $ B$.
Furthermore let $Y_t$ denote the factual ``at-risk'' process,
that is, $Y_t = I( B_{t-} = C_{t-} = 0)$.
We assume that there exist
measurable and bounded functions $\psi^0, \psi^K, \psi^L$ and $\psi^A$ such that
%
\begin{equation}\qquad
E_P \biggl[\int_0^T
h_s \,d\tilde B_s \biggr] = E_P \biggl[
\int_0^T h_s Y_s
\bigl( \psi_s^0 + A \psi_s^A
+ L \psi_s^L + \tilde K_{s-}
\psi_s^K \bigr) \,ds \biggr]
\end{equation}
for every bounded and $\tilde\F_t$-predictable process $h$.

We are now able to
identify the controlled direct effect from surgery.
Note that this is just a slight variation of the model
considered in \cite{MartinussenVansteelandt}.
%
\begin{lemma}
If $\sigma^1$ and $\sigma^2$ are two $\F_0^A \vee\F_t^{\tilde
B}$-predictable processes such that
\begin{eqnarray*}
 E_P \biggl[L \int_0^T
h_t Y_t \exp \biggl( \int_0^t
K_{-s} \psi_s^K \,ds \biggr) \,dt \biggr]& =&
E_P \biggl[ \int_0^T
h_t \sigma^1_t \,dt \biggr],
\\
E_P \biggl[ \int_0^T
h_t Y_t \exp \biggl( \int_0^t
K_{-s} \psi_s^K \,ds \biggr) \,dt \biggr]& =&
E_P \biggl[ \int_0^T
h_t \sigma^2_t \,dt \biggr]
\end{eqnarray*}
for every
$\F_0^{A, B, C}$-predictable and bounded process $h$, then
%
\begin{eqnarray}
&& E_{P_\theta} \biggl[\int_0^T
g_t Y_t \,dB_t \biggr]
\nonumber
\\[-8pt]
\\[-8pt]
\nonumber
&&\qquad = E_{P_\theta}
\biggl[\int_0^T g_t Y_t
\biggl( \psi_t^0 + \psi_t^L
\theta^* \frac{\sigma^1_t}{\sigma^2_t} + \theta^* A \psi_t^A +
\theta^* K_{t-} \psi_t^K \biggr) \,dt \biggr]
\end{eqnarray}
for every
$ \F_t^{B, C}$-predictable and bounded process $g$.
\end{lemma}

\begin{pf*}{Sketch of proof}
By Theorem~\ref{baselinexistthm}, there exist
an $\F_0^A$-measurable random variable $W_0^1$ and an
$\F_0^{A, K, L}$-measurable random variable $W_0^2$ such that
\[
\frac{ dP_\theta|_{\F_0}}{dP|_{\F_0}} = W_0^1 W_0^2
\quad\mbox{and}\quad \frac{ dP_\theta|_{\F_0^{L, A }}}{dP|_{\F_0^{L, A }}} = W_0^1.
\]
If $H_1$ is $\F_0^L$-measurable, $\tilde H_1 := E_P [ H_1  |
\F_0^A ]$ and $H_2$ is $\F_0^A$-measurable, then
%
\begin{equation}
\label{G-ident} E_{P_\theta} [ H_1 H_2 ] =
E_P \bigl[ H_1 H_2 W_0^1
\bigr] = E_P \bigl[ \tilde H_1 H_2
W_0^1 \bigr]= E_{P_\theta} [ \tilde H_1
H_2 ].
\end{equation}

Similarly, let $h$ be a bounded and $\tilde\F_t$-predictable
process, and let $\mu_s^B := \tilde Y_s  ( \psi_s^0 + A \psi_s^A
+ L \psi_s^L + K_{s-}
\psi_s^K  ), $ and note that
\begin{eqnarray*}
E_{P_\theta} \biggl[\int_0^T
h_s \,dB_s \biggr] & = &E_{P} \biggl[\int
_0^T h_s \,dB_s
W_T \biggr]
\\
& =& E_{P} \biggl[\int_0^T
h_s W_{s-} \,dB_s \biggr] +
E_{P} \biggl[\int_0^T
h_s d[ B , W]_s \biggr]
\\
& =& E_{P} \biggl[\int_0^T
h_s W_{s-} \,dB_s \biggr]
\\
& =& E_{P} \biggl[\int_0^T
h_s W_{s-} \mu_s^B \,ds \biggr]
\\
& =& E_{P} \biggl[\int_0^T
h_s \mu_s^B \,ds W_T \biggr]
\qquad\mbox{by \cite{JacodShiryaev}, Proposition~I 3.14}
\\
& =& E_{P_\theta} \biggl[\int_0^T
h_s \mu_s^B \,ds \biggr].
\end{eqnarray*}

One can show that there exists an intermediate
probability measure $\tilde P$ on $\tilde\F_T$
such that:
\begin{longlist}[(1)]
\item[(1)]
\[
P_\theta|_{\tilde\F_T} \ll\tilde P \ll P |_{\tilde
\F_T}.
\]
\item[(2)] For every bounded and Borel-measurable function $h$:
\begin{itemize}
\item$E_{\tilde P}  [ h(A)  ] = h( \theta^* A)$, $\tilde P$-a.s.;
\item$E_{\tilde P}  [ h(L)  | \F_0^{A },
] =
E_{ P}  [ h(L)  | \F_0^{A }  ]$;
\item$E_{\tilde P}  [ h(K_0)  | \F_0^{A, L }
] =
E_{ P}  [ h(K_0)  | \F_0^{A , L }  ]$;
\item
$E_{\tilde P}  [ h(B_0)  | \F_0^{A, L, K }
] =
E_{ P}  [ h(B_0)  | \F_0^{A , L, K }  ]$.
\end{itemize}
\item[(3)]
Whenever $h$ is a bounded and $\tilde\F_t$-predictable process,
then:
\begin{itemize}
\item
\[
E_{\tilde P} \biggl[\int_0^T
h_s \,dB_s \biggr] = E_{\tilde P} \biggl[\int
_0^T h_s \mu_s^B
\,ds \biggr];
\]
\item if
$\mu^K$ and $\mu^C$ are $\tilde\F_t$-predictable processes
such that
\begin{eqnarray*}
E_{ P} \biggl[\int_0^T
h_s \,dK_s \biggr] &=& E_{ P} \biggl[\int
_0^T h_s \mu_s^K
\,dU^K_s \biggr],
\\
E_{ P} \biggl[\int_0^T
h_s \,dC_s \biggr] &=& E_{ P} \biggl[\int
_0^T h_s \mu_s^C
\,dU^C_s \biggr],
\end{eqnarray*}
then
\begin{eqnarray*}
E_{\tilde P} \biggl[\int_0^T
h_s \,dK_s \biggr] & =& E_{\tilde P} \biggl[\int
_0^T h_s \mu_s^K
\,dU_s^K \biggr],
\\
E_{\tilde P} \biggl[\int_0^T
h_s \,dC_s \biggr] & = &E_{\tilde P} \biggl[\int
_0^T h_s \mu_s^C
\,dU_s^C \biggr].
\end{eqnarray*}
\end{itemize}
\end{longlist}
Note that by \cite{JacodMultivariate}, Proposition~4.3, there exists
an $\F_t^{A, L , B }$-adapted $\tilde
P$-martin\-gale $\Xi$ such that
\[
\Xi_T = \frac{dP_\theta|_{ \F_T^{A, L , B, C} }} { d\tilde P|_{
\F_T^{A, L , B, C} } },\qquad  [ B, \Xi ] = 0
\]
and
%
\begin{equation}
\label{deterministicK} Y \Xi_{-} = Y \exp \biggl( - \int
_0^\cdot \theta^* K_r
\psi_r^K \,dr \biggr).
\end{equation}

Bayes's formula with predictable projections shows that
%
\begin{equation}
\label{skewexpectation} E_{\tilde P} \biggl[L \int_0^T
Y_s h_s \,ds \biggr] = E_{\tilde P} \biggl[\int
_0^T Y_s h_s
\frac{ \sigma^1_s }{\sigma^2_s } \,ds \biggr]
\end{equation}
for every bounded and $\F_t^{A, B, C }$-predictable process $h$.
Now,
\begin{eqnarray*}
E_{P_\theta} \biggl[ \int_0^T L
h_s Y_s \,ds \biggr] & =& E_{\tilde P} \biggl[ \int
_0^T L h_s Y_s \,ds
\Xi_T \biggr]
\\
& =& E_{\tilde P} \biggl[ \int_0^T
\Xi_{s-} L h_s Y_s \,ds \biggr]
\\
& = &E_{\tilde P} \biggl[ \int_0^T
\Xi_{s-} h_s Y_s \frac{ \sigma^1_s
}{\sigma^2_s } \,ds
\biggr]\qquad \mbox{by \eqref{deterministicK}}
\\
& = &E_{\tilde P} \biggl[ \int_0^T
h_s Y_s \frac{ \sigma^1_s }{\sigma^2_s } \,ds \Xi_T
\biggr]
\\
& =& E_{ P_\theta} \biggl[ \int_0^T
h_s Y_s\frac{ \sigma^1_s }{\sigma^2_s } \,ds \biggr]
\end{eqnarray*}
for every bounded $\F_t^{A, B, C }$-predictable process $h$,
which implies that \eqref{G-ident} holds.
\end{pf*}

\subsubsection{Consistency of the modified sequential $G$-estimator}

We are now able to show that the modified sequential $G$-estimator
suggested in
\cite{MartinussenVansteelandt} is uniformly consistent, also when we
consider a time-dependent mediating treatment $K$. Let
$\theta_1, \theta_2$ be two actions as in the previous proposition,
but where $\theta_1^* A = 0$ and $\theta^*_2 A = 1$ and consider
corresponding $\F_t^{A, B ,C}$-predictable processes $\gamma^1$ and
$\gamma^2$ as the fractions in \eqref{skewexpectation}.\vadjust{\goodbreak} Furthermore,
we assume that our observations consist of the
event histories for $n$ independent equally distributed
individuals, following the current generic model.
We will also slightly misuse the notation and let $N$, from now on,
denote the corresponding counting process that is aggregated over
the $n$ independent individuals.

\begin{lemma}Let $\widehat\Psi^0 , \widehat\Psi^A, \widehat\Psi^L$
and $\widehat\Psi^K$ denote the usual additive regression
estimators of Aalen, let $\tilde Y := Y_t^B Y_t^C$
and define
\begin{eqnarray*}
\tilde M_t& :=& N_t^B -\int
_0^t \mu_s^B \,ds,
\overline\gamma_t := %
\pmatrix{ \gamma^1_t
\vspace*{2pt}\cr
\gamma^2_t } %
,\qquad
\Gamma_t := %
\pmatrix{\displaystyle \Psi^0_t
+ \int_0^t \rho^1_s
\,d\Psi^L_s
\vspace*{2pt}\cr
\displaystyle\Psi^A_t + \int_0^t
\rho^2_s - \rho^1_s d
\Psi^L_s }, %
\\[-2pt]
\widehat H_t &:=& \operatorname{diag} \pmatrix{\displaystyle \tilde Y^1_t
\exp \biggl( \int_0^t K_{s-}^1
\,d \widehat{\Psi}_s^K \biggr)
\vspace*{2pt}\cr
\vdots\vspace*{2pt}
\cr
\displaystyle\tilde Y^n_t \exp \biggl( \int
_0^t K_{s-}^n \,d
\widehat{\Psi}_s^K \biggr) },\\[-2pt]
 H_t
&:=& \operatorname{diag} \pmatrix{\displaystyle \tilde Y^1_t \exp \biggl( \int
_0^t K_{s-}^1 \,d
\Psi_s^K \biggr) \vspace*{2pt}
\cr
\vdots\vspace*{2pt}
\cr
\displaystyle\tilde Y^n_t \exp \biggl( \int_0^t
K_{s-}^n \,d \Psi_s^K \biggr) },
\\[-2pt]
Z_t& :=& \tilde Y_t \cdot \pmatrix{ 1 & A_1
\vspace*{2pt}
\cr
\vdots& \vspace*{2pt}
\cr
1 & A_n} ,\qquad
Z_{s-}^{\widehat H}:= \bigl( Z_{s-}^T
\widehat H_{s-} Z_{s-} \bigr)^{-1}
Z_{s-}^T \widehat H_{s-},
\\[-2pt]
Z_{s-}^{ H}&:= &\bigl( Z_{s-}^T
H_{s-} Z_{s-} \bigr)^{-1} Z_{s-}^T
H_{s-} \quad\mbox{ and }\\[-2pt]
 \widehat{\Gamma}_t& : =& \int
_0^t Z_{s-}^{\widehat H} \,d
N_s - \int_0^t
Z_{s-}^{\widehat H} K_{s-} \,d \widehat
\Psi_s^K.
\end{eqnarray*}

We have that
%
\begin{equation}
\lim_{n \rightarrow\infty} P \Bigl( \sup_{t \leq T} |\widehat{
\Gamma}_t - \Gamma_t | \geq\delta \Bigr) = 0
\end{equation}
for every $\delta> 0$.
\end{lemma}

\begin{pf}
Note that
%
\begin{eqnarray}
\widehat{\Gamma}_t - \Gamma_t& = & \int
_0^t Z_{s-}^{\widehat H}
\,dN_s - \int_0^t
Z_{s-}^{\widehat H} K_{s-} \,d\widehat{
\Psi}_s^K - \Gamma_t
\\[-2pt]
&= & \label{countleng} \int_0^t
Z_{s-}^{\widehat H} - Z_{s-}^{ H}
\,dN_s + \int_0^t \bigl(
Z_{s-}^{ H} - Z_{s-}^{ \widehat H} \bigr)
K_{s-} \,d \Psi_s^K\vadjust{\goodbreak}
\\
\label{mart} &&{}+ \int_0^t Z_{s-}^{ H}
\,d\tilde M_s + \int_0^t
Z_{s-}^{ \widehat H} K_{s-} \,d \bigl( \Psi^K_s
- \widehat{\Psi}_s^K\bigr)
\\
&&{} \label{loln} + \int_0^t
Z_{s-}^{ H} \pmatrix{ Z_{s-} & L_{s-}
} \,d \pmatrix{\Psi^0_s \vspace*{2pt}
\cr
\Psi^A_s \vspace*{2pt}
\cr
\Psi_s^L
} - \Gamma_t.
\end{eqnarray}

Let
\[
V = \pmatrix{ 1 & 0 \vspace*{2pt}
\cr
-1 & 1}.
\]
We have that $V^T Z^T_{s-} H_{s-} Z_{s-} V = S_{t-}$ where $S_{t-}$ is
a $2\times2$-diagonal matrix. Moreover,
$ ( Z^T_{s-} H_{s-} Z_{s-}  )^{-1} = V S_{t-} V^T$.

Note that $  | \int_0^\cdot Z_{s-}^{ H}
\,d\tilde M_s  |_2^2 $ is Lenglart dominated by $ \operatorname{Tr} \langle
\int_0^\cdot Z_{s-}^{ H}
\,d\tilde M_s \rangle$ and
\begin{eqnarray*}
&& \operatorname{Tr} \biggl\langle \int_0^\cdot
Z_{s-}^{ H} \,d\tilde M_s \biggr
\rangle_T
\\
&& \qquad=  \int_0^T \operatorname{Tr} \bigl(
Z_{s-}^T H_{s-} Z_{s-}
\bigr)^{-1} Z_{s-}^T H_{s-}
\operatorname{diag} \mu H_{s-} Z_{s-} \bigl( Z_{s-}^T
H_{s-} Z_{s-} \bigr)^{-1} \,ds
\\
&&\qquad \leq \int_0^T \operatorname{Tr} \bigl(
Z_{s-}^T H_{s-} Z_{s-}
\bigr)^{-1} \| \operatorname{diag} \mu H_{s-} \|_{\operatorname{op}} \,ds,
\end{eqnarray*}
which converges in probability to $0$.
By Lenglart's inequality \cite{JacodShiryaev}, we obtain that $ \int_0^\cdot Z_{s-}^{ H}
\,d\tilde M_s$ converges uniformly to $0$ in probability with respect to
$P$.

Since
\[
\lim_{\delta\rightarrow\infty} P \Bigl( \sup_s | Z_s
K_s | \geq \delta \Bigr) = 0 \quad\mbox{and}\quad \lim_{\delta\rightarrow\infty} P
\Bigl( \sup_s | \widehat Z_s K_s | \geq
\delta \Bigr) = 0,
\]
and $\widehat
\Psi^K$ converges uniformly in probability to $\Psi^K$, we also
have that
\[
\int_0^t Z_{s-}^{ \widehat H}
K_{s-} \,d \bigl( \Psi^K_s - \widehat
\Psi_s^K\bigr) \quad\mbox{and}\quad \int_0^t
\bigl( Z_{s-}^{ H} - Z_{s-}^{ \widehat H} \bigr)
K_{s-} \,d \Psi^K_s
\]
converge uniformly in probability to $0$ w.r.t. $P$. This shows that
\eqref{mart} converges uniformly in probability to $0$ w.r.t. $P$ as
well.

We have that
\begin{eqnarray*}
Z^H_{s-} L_{s-} = V S_{s-}^{-1}
\pmatrix{\displaystyle \sum_{i = 1}^n
H^i_{s-} L^i_{s-} \bigl( 1-
A^i\bigr) \vspace*{2pt}
\cr
\displaystyle\sum_{i = 1}^n
H^i_{s-} L^i_{s-} A^i }
= V \pmatrix{ \displaystyle\frac{\sum_{i = 1}^n H^i_{s-} L^i_{s-} ( 1- A^i)} {\sum_{i = 1}^n
H^i_{s-} ( 1- A^i)} \vspace*{2pt}
\cr
\displaystyle\frac{\sum_{i = 1}^n H^i_{s-} L^i_{s-} A^i }{\sum_{i = 1}^n H^i_{s-} A^i} }.
\end{eqnarray*}

The law of large numbers implies that $ Z^H_{s-} L_{s-}$ converges
in $P$-probability to $ V
\overline \gamma(s)$.
Now, \eqref{loln} equals
%
\begin{equation}
\int_0^t Z^H_{s-}
L_{s-} - V \gamma(s) \,d \Psi_s^L
\end{equation}
and
%
\begin{eqnarray}
&& E_P \biggl[ \sup_t \biggl| \int_0^t
Z^H_{s-} L_{s-} - V \gamma(s) \,d
\Psi_s^L\biggr | \biggr]
\nonumber
\\[-8pt]
\\[-8pt]
\nonumber
&&\qquad\leq\int_0^T
E_P \bigl[ \bigl| \bigl( Z^H_{s-}
L_{s-} - V \gamma(s-) \bigr) \bigr| \bigr] \bigl| \psi_s^L\bigr|
\,ds.
\end{eqnarray}
Therefore \eqref{loln} converges uniformly in probability w.r.t. $P$.

A computation shows that
$ |\int_0^\cdot Z_{s-}^{\widehat H}
- Z_{s-}^{ H} \,dN_s|_1$ is Lenglart dominated by
\begin{eqnarray*}
&&\|V \|_1 \int_0^\cdot \sum
_j  \biggl| \frac{\widehat H_{s-}^j  ( 1 - A_j
)\mu_s^j }{ \sum_i \widehat H_{s-}^i  ( 1- A_i)
} - \frac{ H_{s-}^j  ( 1 - A_j
)\mu_s^j }{ \sum_i H_{s-}^i  ( 1- A_i)
} \biggr|
\\
&&\hspace*{32pt}\qquad{} +\biggl | \frac{\widehat H_{s-}^j A_j
\mu_s^j }{ \sum_i \widehat H_{s-}^i A_i
} - \frac{ H_{s-}^j A_j\mu_s^j }{ \sum_i H_{s-}^i A_i
} \biggr| \,ds.
\end{eqnarray*}
This process converges uniformly in probability to $0$, so we see that
\eqref{countleng} also converges uniformly in probability to $0$.
This means that $\widehat
\Gamma- \Gamma$ converges uniformly in probability to $0$, so
$\widehat
\Gamma$ actually converges to $\Gamma$ in the similar sense.
\end{pf}

The cumulative $P_{\theta_i}$-hazard of
$\tilde B$ is given by
%
\begin{equation}
\Lambda_t^{\theta_i} = \int_0^t
\psi_s^0 + \theta^*_i A
\psi_s^A + \theta^*_i K_s
\psi_s^K + \gamma_s \psi_s^L
\,ds.
\end{equation}

Since stochastic integrals are continuous with respect to uniform
convergence on compacts in probability, we see that
\begin{eqnarray*}
\lim_{\delta\rightarrow0} P \biggl( \sup_t \biggl| \int_0^t
\pmatrix{ 1, & \theta^*_i A, & \theta^*_i
K_{s-} } \pmatrix{ d\widehat\Gamma_s \vspace*{2pt}
\cr
d
\widehat \Psi_s^K} - \Lambda_t^{\theta_i}
\biggr| \geq\delta \biggr) = 0,
\end{eqnarray*}
that is, we obtain a consistent estimator of $ \Lambda_t^{\theta_i}$.
A consistent estimator for the controlled direct effect of $A$ on $B$
is given by the second component of~$\widehat\Gamma$.
\section{Discussion}

The primary concern in this paper is the possibility of estimating
parameters for the counterfactual situation from the observational
data, given that the counterfactual model is correct. This comes
mainly down to whether the counterfactual probability is absolutely
continuous with respect to the factual probability and whether the
counterfactual parameters of interest are identifiable.
The previously mentioned related works by Arjas and Parner, \cite{Arjas3}
and \cite{ArjasParner}, construct counterfactual probability
distributions by piecing stochastic intervals together as in
\cite{JacodMultivariate}, Section~3. Unlike Parner and Arjas, we
take a more martingale oriented approach, also based on the
seminal paper \cite{JacodMultivariate}. This enables us to apply
directly already well-established methods from stochastic analysis and
martingale theory.
In fact, surprisingly much causal inference
can be well understood in terms of martingale measures, Bayes's
rule and Girsanov's theorem.
This approach translates directly the problem about data re-weighting
into a
thoroughly studied problem in the literature, that is, whether the
stochastic exponential of a local martingale defines a martingale, see
\cite{LepingleMemin} and \cite{ShiryaevKallsen}.

Another difference from the work of Arjas and Parner is that we
consider an explicit intervention in terms of a transformation
$\theta$ on sample space. While not being absolutely necessary, it
still provides additional clarification, as it makes the notion
of counterfactual outcomes more explicit, or
perhaps even demystified. The notation $\operatorname{do}(X=x)$, \cite
{Pearl}, is simply interpreted
as the measurable transformation on the sample space that forces
every outcome of $X$ into $x$ and leaves the remaining observations
unchanged. When the action becomes more complex than
just forcing a variable into a fixed value, this interpretation
becomes even more appealing.

The introduction of the transformation $\theta$
sheds some light on another aspect:
One may in fact
think of a causal inference problem as a stochastic control problem,
or
a decision problem, where the
assumptions about the model are kept as modest as possible.
The main objective in stochastic control theory is to find an
optimal intervention strategy and compute the corresponding expected
payoff. Causal inference appears
as a special case of this, in the sense that there one mostly considers
only one intervention strategy, namely the transformation $\theta$,
and aims to compute the expected payoff.

One is often confronted with latent
factors in epidemiological settings. This lack of information
typically yields nonidentifiable effects.
In special situations, one can use graphical arguments to
ensure identifiability of counterfactual
parameters and also provide exact formulas for these. Such examples are
the \textit{back-door formula}, \textit{front-door formula} and
\textit{sequential back-door formula} \cite{Pearl}, Section~3.3.1, 3.3.2, 4.4.3
and \cite{DidelezEichler}.
We show
that we may take advantage of the local independence graphs
to identify causal effects in
event-history analysis.

When the counterfactual effect is possibly
unidentifiable, one may try to compute upper and lower bounds for
this.
This can also be thought of as a control problem where
``the nature'' is allowed to control the latent factors in order to
maximize or
minimize counterfactual effects. This corresponds to an optimization
problem under constraints. The latent variables may only be altered
in such a way that the observable factors maintain the same joint
distribution and also such that some given directed graph constantly
defines a local independence graph. Let $ \mathcal S $ denote the set
of counterfactual distributions
corresponding to these constraints.
The ``causal effect'' would then be sandwiched by
$
\inf_{P' \in\mathcal S}
E_{ P'}  [ \eta] \leq
E_{ P_\theta}  [ \eta]
\leq \sup_{P'' \in\mathcal S}E_{ P''}  [ \eta]
$.

The set $\mathcal S$ may have a somewhat complicated geometry.
If one instead considers the convex hull,
we obtain other, not
necessarily, tight bounds.
\[
\inf_{P' \in\mathrm{conv}( \mathcal S ) } E_{ P'} [ \eta] \leq E_{ P_\theta} [ \eta]
\leq \sup_{P'' \in \mathrm{conv} ( \mathcal S ) }E_{ P''} [ \eta].
\]
These bounds may be computed by allready developed linear programing
techniques. This approach was for instance taken in
\cite{Balke}, but is likely to generalize to more
complicated continuous-time scenarios as well.

\begin{appendix}
\section*{Appendix}\label{app}

\subsection*{Uniqueness of counterfactual distributions}

\renewcommand{\thelemma}{\Alph{section}.\arabic{lemma}}
\setcounter{lemma}{0}
\begin{lemma} \label{uniquenessatbaseline}
There exists at most one counterfactual distribution $P_\theta$ on $\F_0$ that
imposes contemporaneously independent outcomes.
\end{lemma}

\begin{pf}
Let $T_1, \ldots, T_m$ be an enumeration of $\{T(V)\}_{V
\in\V}$ such that $j < k$ implies $T_j < T_k$.

Assume that $P'$ and $P''$ are two counterfactual distributions that have
contemporaneously independent outcomes and $\eta$ is an
$\F_0^{V_k}$-measurable random variable.
Let $\{X_i\}_i$ be an enumeration of $\{V \in\X| T(V) =
T_1\}$ and let $\{A_j\}_j$ be an enumeration of $\{V \in\A| T(V) =
T_1\}$. Whenever $\{h_i\}_i$ and $\{g_j\}_j$ are two families of
bounded and measurable functions, then
\begin{eqnarray*}
 E_{P'} \biggl[\prod_i
h_i( X_i) \prod_j
g_l( A_j) \biggr]& =& E_{P'} \biggl[\prod
_i h_i( X_i) \biggr]
E_{P'} \biggl[ \prod_j g_j(
A_j) \biggr]
\\
&= & \prod_i E_{P'} \bigl[
h_i( X_i) \bigr] E_{P'} \biggl[ \prod
_j g_j( A_j) \biggr]\\
&=& \prod
_i E_{P''} \bigl[ h_i(
X_i) \bigr] E_{P''} \biggl[ \prod
_j g_j( A_j) \biggr]
\\
&= & E_{P''} \biggl[\prod_i
h_i( X_i) \biggr] E_{P''} \biggl[ \prod
_j g_j( A_j) \biggr] \\
&=&
E_{P''} \biggl[\prod_i h_i(
X_i) \prod_j g_j(
A_j) \biggr].
\end{eqnarray*}
This shows that if $\eta$ is a bounded random variable that only
depends on the information at $T_1$, then $E_{P'} [ \eta] = E_{P''} [
\eta]$.
We continue with an induction argument and assume that
$
E_{P'}  [ \eta ] = E_{P''}  [ \eta ]
$
for every bounded and random variable $\eta$ that only depends on $\{ V
\in\V|T(V) < T_k \}$ and aim to prove that this also holds if
$\eta$ depends on the information at time $T_k$. Let
$\{X_i\}_i$ be an enumeration of $\{V \in\X| T(V) =
T_k\}$, and let $\{A_j\}_j$ be an enumeration of $\{V \in\A| T(V) =
T_k\}$. Whenever $\{h_i\}_i$ and $\{g_j\}_j$ are two families of
bounded and measurable functions, then
\begin{eqnarray*}
&& E_{P'} \biggl[\eta\prod_i
h_i( X_i) \prod_j
g_j( A_j) \biggr]
\\
&&\qquad=  E_{P'} \biggl[\eta E_{P'} \biggl[ \prod
_i h_i( X_i) \Big| \F_0^{p(V_1)}
\biggr] \prod_j \theta^* g_j(
A_j) \biggr]
\\
&&\qquad=  E_{P'} \biggl[\eta \prod_i
E_{P'} \bigl[ h_i( X_i) |
\F_0^{p(V_1)} \bigr] \prod_j
\theta^* g_j( A_j) \biggr]
\\
&&\qquad=  E_{P''} \biggl[\eta \prod_i
E_{P''} \bigl[ h_i( X_i) |
\F_0^{p(V_1)} \bigr] \prod_j
\theta^* g_j( A_j) \biggr]
\\
&&\qquad=  E_{P''} \biggl[\eta E_{P''} \biggl[ \prod
_i h_i( X_i) \Big| \F_0^{p(V_1)}
\biggr] \prod_j \theta^* g_j(
A_j) \biggr]
\\
&&\qquad=  E_{P''} \biggl[\eta \prod_i
h_i( X_i) \prod_j
g_j( A_j) \biggr].
\end{eqnarray*}
This proves the induction hypothesis, that is, $E_{P'} [ \eta] = E_{P''}
[ \eta]$ whenever $\eta$ depends on $\{V \in\A| T(V) \leq
T_k\}$.
\end{pf}
%

\begin{theorem} \label{uniquenessduringfollow-up}
There exists at most one probability measure on $\F_T$ that
simultaneously satisfies
\eqref{baselineaction2}, \eqref{baselineaction1}, \eqref
{follow-upaction1} and \eqref{follow-upaction2}.
\end{theorem}

\begin{pf}
Recall definition \eqref{nutheta}.
\eqref{follow-upaction1} and \eqref{follow-upaction2} imply that
%
\begin{equation}
\label{follow-upaction3} E_{P_\theta} \biggl[\int_J
\int_0^T h(s , x) N( ds, dx) \biggr] =
E_{P_\theta} \biggl[\int_J \int_0^T
h(s , x) \nu^\theta( ds, dx) \biggr].
\end{equation}
Now \cite{JacodMultivariate}, Theorem~3.4, implies that there
exists at most one probability measure on $\F_T$ that coincides with
$P_\theta$ on $\F_0$ and satisfies \eqref{follow-upaction3}.
\end{pf}

\subsection*{Dual predictable projections}

\begin{lemma} \label{predictableprojections}
Let $U$ denote the dual predictable projection of $N$ with respect
to $Q$ onto the filtration $\F_t$.
\begin{longlist}[(1)]
\item[(1)] If $h$ is a bounded and $\PP^V$
measurable processes, then
\[
\int_{J_V} \int_0^\cdot h(
s, x) U(ds,dx)
\]
defines an $\F_t^V$-predictable process of finite variation.
\item[(2)] If $h$ and $h'$ are bounded and $\PP\otimes\J$
measurable processes, then
%
\begin{eqnarray}
 \label{strongort} \biggl[ \int_{J_V} \int
_0^\cdot h( s, x) U(ds,dx), \int
_{J_{V'}} \int_0^\cdot
h'( s, x) U(ds,dx) \biggr] &=& 0,
\\
\label{strongort2}
 \biggl[ \int_{J_V} \int
_0^\cdot h( s, x) U(ds,dx), \int
_{J_{V'}} \int_0^\cdot
h'( s, x) N(ds,dx) \biggr] &=& 0
\end{eqnarray}

$Q$-a.s. whenever $V \neq V'$.
\item[(3)] There exists a nonnegative and $\PP\otimes\J$-measurable
process $\lambda$ such that
\[
E_P \biggl[ \int_J \int_0^T
h( s, x) N(ds, dx) \biggr] = E_P \biggl[ \int_J
\int_0^T h( s, x)\lambda(s, x) U(ds, dx)
\biggr]
\]
for every bounded and $\PP\otimes\J$-measurable
process $h$.
\end{longlist}
\end{lemma}
\begin{pf}
The integral equation
%
\begin{equation}
\int_J \int_0^T h( s ,
x) N^V ( ds, dx) = \int_{J_V} \int
_0^T h( s , x) N ( ds, dx)
\end{equation}
defines a multivariate point process $N^V$ with mark space $J$ which
only jumps at marks in $J_V$. \cite{JacodMultivariate}, Theorem~2.1, provides a dual predictable
projection $U^V$ of $N^V$ with respect to the reference measure $Q$
onto the filtration $\F_t^V$.

Let $h$ be a bounded and $\PP\otimes\J$ measurable process.
\cite{JacodShiryaev}, Theorem~I 2.2.ii and a monotone class
argument provides a bounded and $\PP^V$-measurable process
$h^V$ such that
\[
\tilde h(\cdot, \cdot) = E_Q \bigl[ h( \cdot, \cdot) |
\F_T^V \bigr], \qquad Q\mbox{-a.s.}
\]
Now,
\begin{eqnarray*}
E_Q \biggl[\int_{J_V} \int_0^T
h(s, x) U( ds, dx) \biggr] & =& E_Q \biggl[\int_{J_V}
\int_0^T h(s, x) N( ds, dx) \biggr]
\\
& = &E_Q \biggl[\int_J \int
_0^T h(s, x) N^V( ds, dx) \biggr]
\\
& =& E_Q \biggl[\int_J \int
_0^T \tilde h(s, x) N^V( ds, dx)
\biggr]
\\
&= &E_Q \biggl[\int_J \int
_0^T \tilde h(s, x) U^V( ds, dx)
\biggr]
\\
&= & E_Q \biggl[\int_J \int
_0^T h(s, x) U^V( ds, dx) \biggr],
\end{eqnarray*}
which proves the first claim.

To prove \eqref{strongort}, let $W \subset J_V$ and $W' \subset
J_{V'}$ be measurable subsets and consider the corresponding
counting processes
\[
N^W_t := N\bigl([0, t], W\bigr) \quad\mbox{and}\quad
N^{W'}_t := N\bigl([0,t], W'\bigr)
\]
and let
\[
U^W_t := U\bigl([0, t], W\bigr)\quad \mbox{and}\quad
U^{W'}_t := U\bigl([0,t], W'\bigr).
\]
Following \cite{JacodMultivariate}, Proposition~2.3, we see
that
\[
\Delta U^W_s = E_Q \bigl[ \Delta
N^W_s | \F_{s-} \bigr] \quad\mbox{and}\quad \Delta
U^{W'}_s = E_Q \bigl[ \Delta
N^{W'}_s | \F_{s-} \bigr],\qquad  Q\mbox{-a.s.}
\]

Now,
\begin{eqnarray*}
0 \leq E_Q \bigl[ \bigl[ U^W , U^{W'}
\bigr]_T \bigr] & =& E_Q \biggl[ \sum
_{s \leq T} \Delta U^W_s \Delta
U^{W'}_s \biggr]
\\
&\leq& \sum_{s \leq T} E_Q \bigl[ \Delta
U^W_s \Delta U^{W'}_s \bigr]
\qquad\mbox{by Fatou's lemma}
\\
& =& \sum_{s \leq T} E_Q \bigl[
E_Q \bigl[ \Delta N^W_s |
\F_{s-} \bigr] E_Q \bigl[ \Delta N^{W'}_s
| \F_{s-} \bigr] \bigr]
\\
& =& \sum_{s \leq T} E_Q \bigl[ \Delta
N^W_s \Delta N^{W'}_s \bigr]
\\
& =& 0,
\end{eqnarray*}
so $[U^W, U^{W'}] = 0$, $Q$-a.s.

Whenever $f$ and $f'$ are bounded and $\F_t$-predictable processes,
we have%
\begin{equation}
\qquad\biggl[ \int_0^\cdot f_s
\,dU^W_s, \int_0^\cdot
f_s' \,dU_s^{W'} \biggr] = \int
_0^\cdot f_s f_s'\,
d \bigl[ U^W, U^{W'} \bigr]_s = 0,\qquad Q\mbox{-a.s.}
\end{equation}
Equation \eqref{strongort} is therefore satisfied in the special
case with $h = f \cdot\chi_W $ and
$h' = f' \cdot\chi_ { W'} $. The general case now follows from an
application of the Monotone class theorem. Equation
\eqref{strongort2} follows from an almost similar argument.

For the last claim, let $\nu$ denote the dual predictable
projection of $N$ with respect to $P$ onto the filtration $\F_t$
and note that $\nu\ll U$ since $P \ll Q$. The existence of
$\lambda$ then follows directly from
\cite{JacodMultivariate}, Theorem~4.1.
\end{pf}

\end{appendix}

\section*{Acknowledgments}
The author would like to thank Prof. Odd O. Aalen and Prof. Torben
Martinussen for very helpful discussions on this project.

%


\printaddresses

\end{document}